\documentclass[10pt,reqno]{amsart}

\usepackage{hyperref}
\usepackage{graphicx}
\usepackage{xcolor}
\usepackage{multimedia}
\usepackage{mathrsfs}

\theoremstyle{theorem}

\theoremstyle{theorem}

\theoremstyle{theorem}
\newtheorem{coroll}{Corollary}[section]
\theoremstyle{theorem}
\newtheorem{proposition}{Proposition}[section]
\theoremstyle{definition}
\newtheorem{definition}{Definition}[section]
\theoremstyle{theorem}
\newtheorem{example}{Example}[section]
\theoremstyle{definition}
\newtheorem{remark}{Remark}[section]

\numberwithin{equation}{section}

\newcommand{\EE}{\mathbb{E}}
\newcommand{\RR}{\mathbb{R}}

\newcommand{\FF}{\mathbb{F}}

\usepackage{relsize}
\newcommand{\T}{{\mathsmaller {\rm T}}}

\makeatletter
\newcommand*\bigcdot{\mathpalette\bigcdot@{.5}}
\newcommand*\bigcdot@[2]{\mathbin{\vcenter{\hbox{\scalebox{#2}{$\m@th#1\bullet$}}}}}
\makeatother

\begin{document}
	
	\title[]{On the role of parametrization in models with a misspecified nuisance component}

	\author{H.~S.~Battey}
	\thanks{Department of Mathematics, Imperial College London, 180 Queen's Gate, London SW7 2AZ, UK.}
	\address{Department of Mathematics, Imperial College London.}
	\email{h.battey@imperial.ac.uk}
	
	\author{N.~Reid}
	\thanks{Department of Statistical Sciences, University of Toronto, 700 University Ave, 9th Floor, Toronto, Ontario M5G 1Z5, Canada.}
	\address{Department of Statistical Sciences, University of Toronto.}
	\email{nancym.reid@utoronto.ca}

	\maketitle
	
	\begin{abstract}
		The paper is concerned with inference for a parameter of interest in  models that share a common interpretation for that parameter but that may differ appreciably in other respects. We study the general structure of models under which the maximum likelihood estimator of the parameter of interest is consistent under arbitrary misspecification of the nuisance part of the model. A specialization of the general results to matched-comparison and two-groups problems gives a more explicit and easily checkable condition in terms of a new notion of symmetric parametrization, leading to a broadening and unification of existing results in those problems.  The role of a generalized definition of parameter orthogonality is highlighted, as well as connections to Neyman orthogonality. The issues involved in obtaining inferential guarantees beyond consistency are briefly discussed. 
		
				\medskip \emph{Some key words:} causality; model formulation; model structure; parameter orthogonalization; symmetric parametrization; treatment effect.
	\end{abstract}

	\bigskip

	\subsection*{Significance statement}
	
		Statistical models are often chosen based on a combination of scientific understanding and flexibility or mathematical convenience. While the aspects of core scientific relevance may be relatively securely specified in terms of interpretable interest parameters, the rest of the formulation is often chosen somewhat arbitrarily. In many statistical models, this formulation includes so-called nuisance parameters, which are of no direct subject-matter concern, but are needed to complete the model or reflect the complexity of the data. This paper contributes to the foundations of statistics by studying the interplay between model structure and likelihood inference, under misspecification of the nuisance component.

	\bigskip
	
\section*{Introduction}\label{secIntro}

	Scientific practice in nearly every field relies on the use of mathematical models: a provisional base for describing how observable quantities materialize. Some models, such as Einstein's theory of gravitation, provide predictions that can be verified with such accuracy that they are often considered to be truth. This is, however, an atypical situation: the processes underlying most observable data are much too complicated to be described by exact laws and often require a probablistic element. A useful statistical model captures the essence of the data-generating mechanism in a way that can accommodate new data collected under slightly different conditions. When the primary purpose of a model is to provide insight and explanation, as is typical in many scientific contexts, stable aspects should ideally be capable of interpretation. For this, statistical models will usually have a relatively small number of scientifically interpretable parameters of interest. So-called nuisance parameters are used to complete the specification of the model. 
	
	In contexts where prediction is the main focus, such as arise in many applications of machine learning, model-free approaches are sometimes advocated. This is justified on the grounds that faithful representation of a data-generating process is both unnecessary for such a task and also unrealistically optimistic with very large or complex sets of data. For instance in commercial recommender systems the criteria of speed and accuracy of predictions are key, and the identification of interpretable parameters unimportant. In sharp contrast are contexts in which  the goal is to understand how a treatment or exposure influences the distribution of outcomes, having accounted, to the extent feasible, for the diversity of individuals. Breiman \cite{Breiman2001} described this separation succinctly as the ``two cultures'' of statistical thought. The development in this paper is closer to the view expressed by Cox \cite{Cox2001} in the discussion of that work. 
			
	We consider sets of statistical models that share a common interpretation for the parameter of subject-matter interest but differ in their specification of the nuisance component. A fundamental question is whether reliable inference for the interest parameter is still achievable using standard approaches when the nuisance aspect is misspecified. An exemplar setting is in understanding the effect of treatment or exposure on the health outcomes of individuals. These individuals may exhibit considerable heterogeneity and complex interdependencies stemming, for instance, from certain shared genetic traits. Modelling such interdependencies is a formidable challenge and often tackled by introducing a parameter for each individual, postulated to have been drawn from a parametric distribution. Such doubly-stochastic models are sometimes called random-effects models. Since the choice of parametric distribution for the individual-specific nuisance parameters is typically based on mathematical convenience, the prospect of misspecification seems likely.

	The extensive literature on inference under model misspecification appears to have started with Cox \cite{Cox1961}, who considered two non-nested models for a vector random variable $Y$, specifying joint density or mass functions $m$ and $\check{m}$, say, where the latter is completely known up to a finite-dimensional parameter $\theta$.  
	With $\hat{\theta}$ the  maximum likelihood estimator under model $\check{m}$, Cox \cite{Cox1961} showed that when the true density function is $m$, $\hat{\theta}$ is asymptotically normal with expected value $\theta_m^0$ and covariance matrix given by what is now known as the sandwich formula. The quantity $\theta_m^0$ solves
	\[
	\mathbb{E}_m[\nabla_\theta \log\check{m}(Y;\theta)]_{\theta=\theta_m^0} = 0,
	\]
	where $\mathbb{E}_m$ is expectation under model $m$. Equivalently $\theta_m^0$ minimizes with respect to $\theta$ the  Kullback-Leibler divergence
	\[
	\int m(y)\log\biggl(\frac{m(y)}{\check{m}(y;\theta)}\biggr)dy.
	\]
	More rigorous discussions of the distribution theory and regularity conditions were provided by Huber \cite{Huber1967}, Kent \cite{Kent1982}, and White \cite{White1982a}, \cite{White1982b}. 
	
	Inferential guarantees for the quantity $\theta_m^0$ follow directly under classical regularity conditions on the family of models. Our interest, by contrast, is in estimation of the true value of an interest parameter, which is assumed to have a common interpretation in both the true distribution and in the fitted model. Our focus is on uncovering structure for which standard likelihood theory retains its first-order theoretical guarantees, providing foundational insight into the limits of likelihood inference. Reference to the unsolved nature of this question is widespread; see, for instance, Evans \& Didelez \cite[\S 5]{ED2023}. 	
	
	A different line of exploration, with connections to the double-robustness literature (e.g.~Robins et al.~\cite{Robins94},  and Chernozhukov et al.~\cite{Cherno2018}) is the so-called assumption-lean approach to modelling of Vansteelandt and Dukes~\cite{VD2022}, who centre their analyses on model-free estimands rather than parameters of a given model. The ideal estimand coincides with a particular model parameter when the chosen model includes $m$, and retains a degree of interpretability when the model is misspecified. The tension between estimands and model parameters has a long history, going back to Fisher and Neyman; see Cox \cite{Cox2012}.
	
	In the discussion of \cite{VD2022}, Battey \cite{B2022} conjectured a condition under which the maximum likelihood estimator for an interest parameter is consistent in spite of arbitrary misspecification of the nuisance part of the model. The intuition for that incomplete claim came from a highly involved calculation in a particular paired-comparisons model studied in Battey and Cox~\cite{BC2020}. Any role played by the assumed model for the random effects was unclear from that calculation, and is clarified by the more general insight provided here. 
	
	In \S \ref{secConsistencyGeneral} we study the general structure for such a consistency result, and specialize the analysis to settings of common relevance, providing more explicit conditions on the parametrization for matched-pair and two-group problems. An important conclusion from this analysis is that the conditions for consistency can sometimes be checked without knowledge of the true model. Section \ref{secExamples} shows through this route that some key insights have been overlooked in a rather large literature on misspecified random effects distributions, particularly the role of parametrization. Subsection \ref{secExpFamily} recovers two well-known results for generalized linear models as an illustration of the more general statements, and \S \ref{secMSM} provides a recent example in the context of a marginal structural model, for which the more general results of \S \ref{secConsistencyGeneral} provide elucidation. While some of the material in \S \ref{secConsistencyGeneral} is rather technical, we have aimed to provide intuition for the main results and have provided several examples in \S \ref{secExamples}. The examples are deliberately detached from specific subject-matter details, in order that they be broadly useful for a range of scientific applications.
	
	\section{Consistency in general misspecified models}\label{secConsistencyGeneral}
	
	\subsection{General conditions for consistency}
	
	Let $m$ represent the density function for the outcomes, parametrized in terms of an interest parameter $\psi$ with true value $\psi^*$. The assumed model, while sharing the same interpretable interest parameter, is misspecified in other ways. The joint density function for the observations under the assumed model has parameters $(\psi,\lambda)\in\Psi\times \Lambda$, and $\ell(\psi,\lambda)= \log \check{m}(y;\psi,\lambda)$ denotes the observed log-likelihood function for that model, viewed as a function of $(\psi,\lambda)$ for observed data $y=(y_1,\ldots,y_n)$. As above, maximization of $\ell(\psi,\lambda)$ gives estimates $(\hat{\psi},\hat{\lambda})$; their dependence on $y$ is suppressed in the notation. As functions of the random variable $Y = (Y_1, \dots, Y_n)$ the probability limit of the maximum likelihood estimator, as $n\rightarrow\infty$, is $(\psi_m^0, \lambda_m^0)$, the solution to
	\begin{equation}\label{eqLimit}
	\mathbb{E}_{m}[\nabla_{(\psi,\lambda)}\ell(\psi_m^0, \lambda_m^0)]=0.
	\end{equation}
	We consider the model to be misspecified if there are no values of $\lambda$ in $\Lambda$ for which the true density $m$ is recovered. This precludes the situation in which the true distribution belongs to a submodel of an assumed encompassing family, briefly discussed in \S \ref{secOverstratification}. 
	
	The following example will be used repeatedly to illustrate key ideas and notation.
	
	\begin{example}\label{exRunningExample}
		Consider a matched comparison problem in which, for each of $n$ twin pairs, one individual from each pair is chosen at random to receive a treatment, the other being the untreated control. Let $Y_{i1}$ and $Y_{i0}$ denote the time until some medical event of interest $($e.g.~recovery from illness$)$ for the treated and untreated individual in the $i$th pair. A simple parametric model specifies the outcomes $Y_{i1}$ and $Y_{i0}$ as exponentially distributed of rates $\gamma_i\psi>0$ and $\gamma_i/\psi>0$ respectively. The pair-specific nuisance parameter $\gamma_i$ captures, for instance, genetic differences between the pairs of twins. The parameter of interest $\psi$ with true value $\psi^*$ quantifies the effect of the treatment: $\psi^2$ being the multiplicative effect of the treatment on the instantaneous probability of recovery relative to baseline. Other parametrizations are possible, but the above symmetric parametrization will be shown to play an important role. Suppose the pair effects are modelled as gamma distributed of shape $\kappa$ and rate $\rho$. Then under the assumed model, the joint density function for the outcomes in any given pair at $(y_1,y_0)$ is
		\begin{equation}\label{eqFitted}
		\frac{\Gamma(\kappa+2)\rho^\kappa}{\Gamma(\kappa)(y_1\psi + y_0/\psi + \rho)^{\kappa+2}}.
		\end{equation}
		The true random effects distribution could have been quite different, so that the model \eqref{eqFitted} is misspecified. Nevertheless, the interpretation of the interest parameter $\psi$ is stable over the different specifications. The notional nuisance parameter is $\lambda=(\kappa, \rho)$. 
	\end{example}
	
	\begin{remark}
				A different analysis is possible in this example, treating  the pair effects $\gamma_1,\ldots,\gamma_n$ as fixed arbitrary constants, which can be eliminated from the analysis by basing likelihood inference on the  distribution of the ratios $Z_i=Y_{i1}/Y_{i0}$; this is discused breifely in \S \ref{secFixedRandom}.
	\end{remark}

	\begin{definition}[parameter $m$-orthogonality]
		Let $\nabla^2_{\psi\lambda}\ell(\psi,\lambda)$ denote the cross-partial derivative of the assumed log-likelihood function. The parameter $\psi$ is said to be $m$-orthogonal to the notional parameter $\lambda$ if $\mathbb{E}_m[\nabla^2_{\psi\lambda}\ell(\psi,\lambda)]=0$. This can hold globally for all $(\psi,\lambda)$ or at particular values, the notation $\Psi \perp_m \Lambda$ indicating global $m$-orthogonality and $\Psi \perp_m \lambda$ and $\psi \perp_m \Lambda$ indicating, respectively, local $m$-orthogonality at $\lambda$ for any $\psi$, and local $m$-orthogonality at $\psi$ for any $\lambda$.
	\end{definition}

To gain some intuition for parameter $m$-orthogonality, consider first the much stronger condition, $\nabla^2_{\psi\lambda}\ell(\psi,\lambda)=0$ for all $\psi$ and $\lambda$. This corresponds to what is called a cut in the parameter space, and implies $\mathbb{E}_m[\nabla^2_{\psi\lambda}\ell(\psi,\lambda)]=0$ for any $m$ so that the true density $m$ trivially plays no role. In this strongest setting, parameter orthogonality is a purely geometric property of the log-likelihood function that holds only in certain parametrizations. It is effectively an absence of torsion, implying that, for any fixed $\lambda$, the corresponding log-likelihood function over $\Psi$ is maximized at the same point. This parameter separation in the likelihood function is, however, limited to a relatively small number of models. The weaker condition $\mathbb{E}_m[\nabla^2_{\psi\lambda}\ell(\psi,\lambda)]=0$ only requires the absence of torsion to hold on average over hypothetical repeated draws from the true distribution. If the model is correctly specified, the usual definition of parameter orthogonality with respect to Fisher information is recovered.

	Propositions \ref{propConsistency} and \ref{propCoxWong} give two alternative general conditions for consistency of $\hat{\psi}$ in spite of arbitrary misspecification of the nuisance part of the model. Their proofs, and those of subsequent results, are in Appendix \ref{secProofs}.

	\begin{proposition}\label{propConsistency}
		Let the observed log-likelihood function for the assumed model be strictly concave as a function of $(\psi, \lambda)$. Then $\psi_m^0 = \psi^*$ if and only if $\mathbb{E}_{m}[\nabla_{\psi}\ell(\psi^*, \lambda_m^0)]=0$. The latter condition is equivalent to $
		\mathbb{E}_{m}[\nabla_{\psi}\ell(\psi^*, \lambda)]=0$ for all $\lambda$ if and only if $\psi^* \perp_{m} \Lambda$.
	\end{proposition}
	
	\begin{remark}
		Local orthogonality at a particular value, $\lambda^\prime$ say, is not sufficient to ensure that
		$\mathbb{E}_{m}[\nabla_{\psi}\ell(\psi^*, \lambda^\prime)]=0$ implies $\mathbb{E}_{m}[\nabla_{\psi}\ell(\psi^*, \lambda_m^0)]=0$, as is clear from the proof of Proposition \ref{propConsistency}.
	\end{remark}

	\begin{remark}
	In Example \ref{exRunningExample} it was shown in \cite{BC2020} that $\psi$ is globally $m$-orthogonal to $\lambda=(\kappa,\rho)$ for any distribution over the random effects, and the consistency of the maximum likelihood estimator of $\psi$ was established through an explicit calculation. An initial motivation for the present work was to understand how sensitive that conclusion was to various aspects of the model formulation.
	\end{remark}

	The first part of Proposition \ref{propConsistency}, i.e.~the requirement of an unbiased estimating equation $\mathbb{E}_{m}[\nabla_{\psi}\ell(\psi^*, \lambda_m^0)]=0$, is almost immediate. Since $\lambda_m^0$ depends on the unknown $m$ and is typically unavailable, verification of the orthogonality condition $\psi^* \perp_{m} \Lambda$ uniformly over $m$, in the event of its validity, appears to be the simplest way of establishing consistency. The simplest special cases in which orthogonality can be checked are those in which there is a parameter cut. Such examples include misspecification of the dispersion component of a generalized linear model, discussed in \S \ref{secExpFamily}, and Gaussian linear mixed models with a misspecified distribution over the random effects. Robustness of inference in the latter case has often been observed empirically, e.g.~Schielzeth, et al.~\cite{Sch} without reference to any underlying structure. A more elaborate example with a parameter cut is the causal model of \cite{ED2023} outlined in \S \ref{secMSM}. Example \ref{exRunningExample} obeys the weaker form $\mathbb{E}_m[\nabla^2_{\psi\lambda}\ell(\psi,\lambda)]=0$
	from Proposition \ref{propConsistency}, as does the unbalanced two-group problem of Example \ref{exCWV2} in \S \ref{secMissRE}, in which the pair-specific parameters are replaced by stratum-specific parameters. The relevant structure underpinning these paired and two-group examples is elucidated in \S \ref{secSymmetry}, where some intuition is provided. 
    
    The unbalanced two-group problem of Example \ref{exCW}, on the other hand, violates the $m$-orthogonality assumption of Proposition \ref{propConsistency} with respect to one of its two nuisance parameters. Nevertheless, Cox and Wong \cite{CoxWong2010} showed that the maximum likelihood estimator of the interest parameter in Example \ref{exCW} is consistent, suggesting that Proposition \ref{propConsistency}, while easier to verify, is too strong for some situations. Proposition \ref{propCoxWong} presents weaker conditions for consistency. 
     
    Suppose that the orthogonality condition $\psi^* \perp_{m} \Lambda$ from Proposition \ref{propConsistency} fails, so that there exists at least one $\lambda\neq \lambda_m^0$ such that $\mathbb{E}_{m}[\nabla_{\psi}\ell(\psi^*, \lambda)]\neq 0$. Let
	\begin{equation}\label{eqg}
	g_{\psi}(\psi^*,\lambda) :=
	\mathbb{E}_{m}[\nabla_{\psi}\ell(\psi^*, \lambda)],  \quad
	g_{\lambda}(\psi^*,\lambda) :=
	\mathbb{E}_{m}[\nabla_{\lambda}\ell(\psi^*, \lambda)],
	\end{equation}
	and partition the inverse of the Fisher information matrix $i:=i(\psi^*,\lambda):=\mathbb{E}_m[-\nabla^2_{(\psi,\lambda)}\ell(\psi^*,\lambda)]$ as
	\begin{equation}\label{eqInv}
	\left(\begin{array}{cc}
	i^{\psi\psi}&  i^{\psi\lambda} \\
	i^{\lambda\psi} & i^{\lambda\lambda}
	\end{array} \right) := 	\left(\begin{array}{cc}
	i_{\psi\psi}&  i_{\psi\lambda} \\
	i_{\lambda\psi} & i_{\lambda\lambda}
	\end{array} \right)^{-1}.
	\end{equation}
	In Proposition \ref{propCoxWong} and its proof, we refer to $g_{\psi}(\psi^*,\lambda)$ and $g_{\lambda}(\psi^*,\lambda)$ in the shorthand $g_\psi$ and $g_\lambda$ respectively.
	
	\begin{proposition}\label{propCoxWong}
		If $i^{\psi\psi}g_\psi +i^{\psi\lambda}g_\lambda = 0$ for all $\lambda\in\Lambda$, then $\psi_m^0=\psi^*$. If $\psi$ and $\lambda$ are both scalar parameters, the condition reduces to $i_{\lambda\lambda} g_\psi = i_{\psi\lambda} g_\lambda$.
	\end{proposition}
	
	That Proposition \ref{propCoxWong} is more general than Proposition \ref{propConsistency} is seen on noting that the orthogonality condition of Proposition \ref{propConsistency} implies $i^{\psi\lambda}=0$ and $g_\psi=0$, so that the condition of Proposition \ref{propCoxWong} also holds. The conditions of Propositions \ref{propConsistency} and \ref{propCoxWong} are clearly not necessary. In particular, if $\lambda_m^0$ can be calculated, the route to establishing consistency or inconsistency of $\hat{\psi}$ is more direct. Any practical relevance, beyond general theoretical understanding, is to situations in which $\lambda_m^0$ is not calculable, typically because $m$ is unknown. The conditions in Propositions \ref{propConsistency} and \ref{propCoxWong}, while dependent on $m$, are nevertheless checkable in some situations, as illustrated in the examples of subsequent sections.

	A generalization involves a second nuisance parameter vector, $\nu\in \mathcal{N}$ say, such that $\psi^* \perp_m \mathcal{N}$ and $\lambda \perp_m \mathcal{N}$ in Proposition \ref{propConsistency} or $\Lambda \perp_m \mathcal{N}$ in Proposition \ref{propCoxWong}. It can be shown via a straightforward extension of the arguments leading to Propositions \ref{propConsistency} and \ref{propCoxWong} that the conclusion of the latter is unchanged by this modification. Examples involving a second nuisance parameter are provided in \S \ref{secExamples}.
	
	The consistency demonstrated by \cite{BC2020} for Example \ref{exRunningExample} arises as a special case of Proposition \ref{propConsistency}, while the argument of Cox and Wong \cite{CoxWong2010} for an unbalanced doubly-stochastic two-group problem turns out to be an application of Proposition \ref{propCoxWong} and Corollary \ref{corollOrthog} below.	We return to this example in \S \ref{secTwoGroup}. 
	
	\subsection{Parameter orthogonality and orthogonalization in misspecified models}\label{secOrthog}
	
	In view of Propositions \ref{propConsistency} and \ref{propCoxWong}, a key question is when the assumed Fisher information matrix, with expectation computed under an erroneous model, coincides with the true Fisher information matrix in some or all regions of the parameter space. Such coincidence implies in particular that parameter orthogonality established under a notional model is true more generally. The practical import is that otherwise $\Psi \perp_m \Lambda$ or its local analogue would typically not be verifiable without knowledge of $m$. A partial answer is presented as Proposition \ref{propOrthog}, whose proof is immediate by the definition of sufficiency. 
	
\begin{proposition}\label{propOrthog}
	Suppose that the cross-partial derivative of the assumed log-likelihood function $\nabla_{\psi\lambda}^2\ell(\psi,\lambda)$ depends on the data only through a sufficient statistic $S=(S_{1},\ldots,S_{k})$, for some $1\le k \le n$. Write $\check{\imath}_{\psi\lambda}(\psi,\lambda)=\mathbb{E}_{(\psi,\lambda)}[-\nabla_{\psi\lambda}^2\ell(\psi,\lambda)]$, where $\mathbb{E}_{(\psi,\lambda)}$ means expectation under the assumed model. Provided that $\nabla_{\psi\lambda}^2\ell(\psi,\lambda)$ is additive in $S_{1},\ldots,S_{k}$ and $\mathbb{E}_m(S_j)=\mathbb{E}_{(\psi,\lambda)}(S_j)$ for all $j$, then $i_{\psi\lambda}(\psi,\lambda)=\check{\imath}_{\psi\lambda}(\psi,\lambda)$. 
\end{proposition}

\begin{coroll}\label{corollOrthog}
	Under the conditions of Proposition \ref{propOrthog}, a sufficient condition for parameter $m$-orthogonality $\Psi \perp_m \Lambda$ is $\check{\imath}_{\psi\lambda}(\psi,\lambda)=0$ for all $(\psi,\lambda)\in\Psi\times\Lambda$, and similarly for local $m$-orthogonality.
\end{coroll}
	
	\begin{remark}
	For Example \ref{exRunningExample} $m$-orthogonality can be established more directly via Proposition \ref{propConsistencyPairs}.
	\end{remark}

	In general, a reparametrization of a model from parameters $\phi$ to $\theta = \theta(\phi)$ is called interest-respecting if the parameter of interest is common to both: if $\phi = (\psi, \xi)$ then $\theta = (\psi, \lambda(\psi,\xi))$. The relevance of interest-respecting reparametrizations in scientific work is that the more important interpretable aspects of the model are retained. In \cite{CoxReid1987} interest-respecting orthogonal reparametrizations were explored for improving inference in correctly specified models. The role in misspecified settings is amplified in view of Proposition \ref{propConsistency}, which concerns the first-order properties of the estimator.
	
	In a slightly more explicit notation, suppose that $\phi$ is a parametrization for which the $(r,s)$th component satisfies $i^{(\phi)}_{rs}(\phi)=\check{\imath}^{(\phi)}_{rs}(\phi)=0$. Consider an interest-respecting reparametrization from an initial parameter $\phi$ to a new parameter $\theta$, with coordinates 
	\[
	\phi^r(\theta^{1},\ldots,\theta^p), \quad \theta^a(\phi^{1},\ldots,\phi^p),   \quad  (a, r=1,\ldots, p),
	\]
	where $\theta^1=\phi^1=\psi$. An implication of Corollary \ref{corollOrthog}  is that the approach of Cox and Reid \cite{CoxReid1987} can still be used in misspecified models, provided that the additional structure of Proposition \ref{propOrthog} is present in the parametrization that is orthogonal under the assumed model. Specifically, from a starting parametrization $\phi$,  an orthogonal parametrization $\theta$ is any solution in $\theta^{2},\ldots,\theta^p$ to the system of $p-1$ differential equations
	\[
	\check{\imath}^{(\phi)}_{ab} \frac{\partial \phi^a(\theta)}{\partial \theta^1}+\check{\imath}^{(\phi)}_{1b} = 0, \quad b=2,\ldots,p,
	\]
	or in matrix notation 
	\[
	\frac{\partial \phi(\theta)}{\partial \theta^1} = - (\check{\imath}^{(\phi)}_{\,\bigcdot \,\bigcdot})^{-1}(\check{\imath}^{(\phi)}_{\,\bigcdot\,1}),
	\]
	where $\check{\imath}^{(\phi)}_{\,\bigcdot \,\bigcdot}$ is the Fisher  information matrix under the assumed model without the row and column corresponding to $\phi^1$, and $\check{\imath}^{(\phi)}_{\,\bigcdot\,1}$ is the excluded column. 
	
	Even under the assumptions of Proposition \ref{propOrthog}, with model misspecification, the asymptotic variance of $\hat{\psi}$ is not  $i^{\psi\psi}$ but is given by the sandwich formula of \cite{Cox1961}. This means that inferential guarantees beyond consistency are not, in general, available using standard likelihood-based approaches; see \S \ref{secInf} for a brief discussion of some exceptions.

	\subsection{Symmetric parametrizations and induced antisymmetry}\label{secSymmetry}
	
	For the matched-pair and two-group examples of \S \ref{secMissRE} and \S \ref{secTwoGroup}, a fundamental role is played by symmetric parametrization in transformation models, which induces an antisymmetry on the log-likelihood derivative and thereby the $m$-orthogonality of Proposition \ref{propConsistency} under arbitrary misspecification of the model. The class of problems covered in this subsection exemplifies settings in which $m$-orthogonality can be straightforwardly verified without knowledge of $m$. A simplified explanation is that in transformation models there is a convenient duality between the parameter space and the sample space which results in cancellation of the relevant terms in $\EE_m[\nabla_{\psi\lambda}^2\ell(\psi,\lambda)]$ provided that the appropriate parameterization is used.

	Let $\mathcal{P}$ be a set of continuous probability measures on $\mathscr{Y}$. A transformation model under the action of $G$ (e.g.~Barndorff-Nielsen and Cox~\cite[p.~53]{B-NC1994}) is a subset of $\mathcal{P}$ parametrized by $\gamma\in \Gamma$ and possibly other parameters suppressed in the notation, such that
	\[
	\mathcal{P}_{G}:=\{p\in\mathcal{P}: p(gE;g\gamma)=p(E;\gamma), \;\; g\in G, \, E\in\mathscr{E}, \, \gamma \in \Gamma\},
	\]
	where $\mathscr{E}$ is the set of measurable events on the sample space $\mathscr{Y}$. The definition implies that $G$ acts on both the sample space and the parameter space. In particular, the action of $G$ on $\FF$, say, is a continuous map $G \times \FF \rightarrow \FF$, $(g,x)\rightarrow gx$ for $x\in\FF$ where, in the transformation model, $x$ represents either the data point $y\in \mathscr{Y}$ or the parameter value $\gamma\in \Gamma$. In general $\mathscr{Y}$ and $\Gamma$ need not be equal, although this will often be the case. The identity element on $\mathscr{Y}$ or on $\Gamma$ is written $e=g^{-1}g=g g^{-1}$, the context leaving no ambiguity. 
	
	The relevant group actions for present purposes depend only on the interest parameter $\psi$, so that we may write, in a more explicit notation, $g=g_\psi\in G$. Two examples serve to illustrate the properties of the transformation models. In location models, $\mathscr{Y}=\Gamma=\RR$, the group action is addition and $g_\psi \gamma = \gamma+\psi$, giving $gg^{-1} = g_{\psi}g_{-\psi}= e = g_0$. In scale models, $\mathscr{Y}=\Gamma=\RR^+$, the group action is multiplication and $g_{\psi}\gamma=\gamma \psi$, giving $gg^{-1} = g_{\psi}g_{(1/\psi)}= e = g_1$. 
	
	Let $U$ be a random variable with a distribution depending only on $\gamma$. This is best thought of as the random variable at baseline with respect to $\psi$, i.e.~$\psi=0$ in location models and $\psi=1$ in scale models. Write the probability density function of $U$ at $u$ as $f_U(u;\gamma)du$. The following definition plays an important role in establishing consistency for a treatment parameter in matched-pair and two-groups problems.

	\begin{definition}[symmetric parametrization]\label{defSymmetric}
		Let $Y_1$ and $Y_0$ be independent random variables with probability measures in $\mathcal{P}_G$ and density functions $f_1$ and $f_0$ respectively. Their joint distribution is said to be parametrized $\psi$-symmetrically with respect to $(\psi,\gamma)$ if $g=g_\psi\in G$ depends only on $\psi$ and if the density functions $f_1$ and $f_0$ relate to $f_U$ by 
		\[
		f_{U}(u;\gamma)du = f_1(gu; g\gamma)d(gu) = f_0(g^{-1}u; g^{-1}\gamma)d(g^{-1}u).
		\]
		In other words, the joint density function, when expressed in terms of $u_1=g^{-1}y_1$ and $u_0=gy_0$, is symmetric in $u_1$ and $u_0$:
		\begin{equation}\label{eqJointSym}
		f_1(y_1; g\gamma)f_0(y_0; g^{-1}\gamma)dy_1dy_0 = f_U(u_1;\gamma)f_U(u_0;\gamma)du_1 du_0.
		\end{equation}
	\end{definition}

	Definition \ref{defSymmetric} says that the random variables $U_1=g^{-1} Y_1$ and $U_0=gY_0$ are equal in distribution to $U$, which has a ``standardized form''. The inverted commas express that $f_U$ might not be a standardized form in the conventional sense, as it depends on $\gamma$ and possibly other parameters suppressed in the notation. The symmetric parametrization can equivalently be written 
			\begin{align}
	\begin{split}\label{eqSymmPairGeneral}
		f_1(y_1; g\gamma)dy_1 =& \,f_U(g^{-1}y_1;\gamma) J_{1}^+ dy_1 \\
		f_0(y_0; g^{-1}\gamma)dy_0 =& \, f_U(gy_0;\gamma)J_0^+ dy_0,
	\end{split}
\end{align}
where 
\[
J_1^+ = |d(g^{-1}y_1)/dy_1|, \quad  J_0^+ = |d(gy_0)/d y_0|
\]
satisfy $J_1^+ J_0^+=1$. In a location model, \eqref{eqSymmPairGeneral} becomes
		\begin{align}
		\begin{split}\label{eqSymmPairLocation}
		f_1(y_1; \gamma+\psi)dy_1 =& f_U(y_1-\psi;\gamma) dy_1 \\
		f_0(y_0; \gamma-\psi)dy_0 =& f_U(y_0+\psi;\gamma) dy_0,
		\end{split}
		\end{align}
		while in a scale model \eqref{eqSymmPairGeneral} becomes
		\begin{align}
		\begin{split}\label{eqSymmPairScale}
		f_1(y_1; \gamma\psi)dy_1 =& f_U(y_1/\psi;\gamma)(1/\psi) dy_1 \\
		f_0(y_0; \gamma/\psi)dy_0 =& f_U(y_0\psi;\gamma)\psi dy_0.
		\end{split}
		\end{align}

	Distributions within the scale family are often parametrized in terms of rate, or inverse scale, which reverses the roles of $g$ and $g^{-1}$ relative to their actions in the scale parametrization. Example \ref{exRunningExample} is of this form.
	
	In a location-scale model the action of $G$ is $g=g_{\psi_2}\circ g_{\psi_1}$ where $g_{\psi_1}$ is multiplication by a positive scalar $\psi_1$ and $g_{\psi_2}$ is addition of a scalar $\psi_2$, its inverse being $g^{-1}=g_{\psi_1}^{-1}\circ g_{\psi_2}^{-1}$. Thus, a two-parameter symmetric parametrization in a location-scale model is in terms of $g\gamma=\psi_1 \gamma + \psi_2$ and $g^{-1}\gamma=(\gamma-\psi_2)/\psi_1$. For the applications we have in mind, the interest parameter $\psi$ represents a treatment effect. Thus, if the treatment is assumed to affect either the location or the scale but not both, a location-scale distribution can effectively be treated as either location or scale. Examples include the extreme-value distributions, the most common parametric families arising in renewal theory and used in survival modelling, and members of the elliptically symmetric family in one or more dimensions, which encompasses the Gaussian, Student-t, Cauchy and logistic distributions. 
	
	Consider, as a function of $\psi$ alone, the log-likelihood contribution of $y_1$ and $y_0$, realizations of $Y_1$ and $Y_0$. This is of the form
	\[
\ell(\psi; \gamma, y_1, y_0) = \log L(\psi; \gamma, y_1, y_0) = \log f_1(y_1; g\gamma) + \log f_0(y_0; g^{-1}\gamma).
\]
	
	Variables to the right of the semi-colon in the log-likelihood function $\ell$ are treated as fixed (although arbitrary), together with any other parameters suppressed in the notation.
	
	\begin{definition}[antisymmetry]\label{defAntisymm}
		The symmetric parametrization of $f_1$ and $f_0$ is said to induce antisymmetry on the associated log-likelihood derivative with respect to $\psi$ if $\nabla_\psi \ell(\psi;\gamma,y_1,y_0)$, when expressed in terms of $u_1=g^{-1}y_1$ and $u_0=g y_0$, satisfies $\nabla_\psi \ell(\psi;\gamma,u_1,u_0)=-\nabla_\psi \ell(\psi;\gamma,u_0,u_1)$. 
	\end{definition}

\begin{example}[Continuation of Example \ref{exRunningExample}]\label{exRunningEx2}
In the exponential matched pair problem, the gamma distribution over the random effects is irrelevant for illustrating the symmetric parametrization and induced antisymmetry, as these definitions are conditional on $\gamma$ within a single pair. Since multiplication of rates corresponds to division of scale, it is natural, for consistency with \eqref{eqSymmPairScale}, to define the group operation $g$ as division by $\psi$ rather than multiplication by $\psi$. Thus let $u_1=g^{-1}y_1 = \psi y_1$ and $u_0=g y_0 =y_0/\psi$. Conditionally on $\gamma$, the joint density function of $(Y_1,Y_0)$ is
\[ 
f_1(y_1\hspace{0.5pt};\psi^*,\gamma)f_0(y_0 \hspace{0.5pt};\psi^*,\gamma) dy_1 dy_0
= \gamma^2 \exp\{-\gamma(u_1+u_0)\} du_1 du_0
\]
which is symmetric in $u_1$ and $u_0$. The log-likelihood derivative with respect to $\psi$ is
\[
\nabla_{\psi} \ell (\psi; u_0,u_1) = -\gamma(y_1 - y_0/\psi) = - \gamma(u_1/\psi - u_0/\psi), 
\]
where the right hand side is $- (- \gamma(u_0/\psi - u_1/\psi))=-\nabla_{\psi} \ell (\psi; u_1,u_0)$. This shows that the log-likelihood derivative is antisymmetric in the sense of Definition \ref{defAntisymm}.
\end{example}

The main transformation models arising in statistics are the location models and the scale models. For these, we show in Example \ref{exLocScale} that the the symmetric parametrization automatically induces antisymmetry on the associated log-likelihood derivative. We also show this for a rotation family on the circle in Example \ref{exCircle2}. We conjecture that antisymmetry according to Definition \ref{defAntisymm} is a necessary consequence of the $\psi$-symmetric parameterization for any transformation model, but in the absence of general group-theoretic proof, we present the following necessary and sufficient conditions for asymmetry, that may be checked on a case-by-case basis for more exotic groups. Examples \ref{exLocScale} and \ref{exCircle2} illustrate the application of Proposition \ref{propGroup}.
	
	\begin{proposition}\label{propGroup}
		Suppose that the joint distribution of $Y_1$ and $Y_0$ is parametrized $\psi$-symmetrically with respect to $(\psi,\gamma)$ in the sense of Definition \ref{defSymmetric}. The parametrization induces antisymmetry on the log-likelihood derivative in the sense of Definition \ref{defAntisymm} if and only if 
		\[
		a(u_1,u_0):=\frac{\partial u_1}{\partial \psi} + \frac{\partial u_0}{\partial \psi} = -a(u_0,u_1)
		\]
		and 
		\[
c(u_1,u_0):=\Bigl(\frac{\partial}{\partial \psi}\Bigl|\frac{\partial u_1}{\partial y_1}\Bigr| \Bigr)\Bigl|\frac{\partial u_0}{\partial y_0}\Bigr| + \Bigl(\frac{\partial}{\partial \psi}\Bigl|\frac{\partial u_0}{\partial y_0}\Bigr| \Bigr)\Bigl|\frac{\partial u_1}{\partial y_1}\Bigr| = -c(u_0,u_1),
\]
		where $u_1=g^{-1}y_1$ and $u_0=g y_0$.
	\end{proposition}
	
	A clearer but more cumbersome expression of the conditions $a(u_1,u_0)=-a(u_0,u_1)$ and $c(u_1,u_0)=-c(u_0,u_1)$ makes the dependence of the partial derivatives on $u_1$ and $u_0$ explicit. In this more explicit notation, $a(u_1,u_0)=-a(u_0,u_1)$ amounts to $(\partial u_1/\partial \psi)(u_1) = -(\partial u_0/\partial \psi)(u_1)$ and vice versa.
	
	\begin{example}\label{exLocScale}
		In a symmetric parametrization of a location model $\partial u_1/\partial \psi = -1= -\partial u_0/\partial \psi$ and 
		\[
		\Bigl(\frac{\partial}{\partial \psi}\Bigl|\frac{\partial u_1}{\partial y_1}\Bigr| \Bigr)\Bigl|\frac{\partial u_0}{\partial y_0}\Bigr| = 0 = - \Bigl(\frac{\partial}{\partial \psi}\Bigl|\frac{\partial u_0}{\partial y_0}\Bigr| \Bigr)\Bigl|\frac{\partial u_1}{\partial y_1}\Bigr|.
		\]
		Thus, ${a(u_1,u_0)=0=-a(u_0,u_1)}$ and $c(u_1,u_0)=0=-c(u_0,u_1)$. In a symmetric parametrization of a scale model ${\partial u_1/\partial \psi = -u_1/\psi}$, ${\partial u_0/\partial \psi = u_0/\psi}$, and
		\[
		\Bigl(\frac{\partial}{\partial \psi}\Bigl|\frac{\partial u_1}{\partial y_1}\Bigr| \Bigr)\Bigl|\frac{\partial u_0}{\partial y_0}\Bigr| = -\frac{1}{\psi} = - \Bigl(\frac{\partial}{\partial \psi}\Bigl|\frac{\partial u_0}{\partial y_0}\Bigr| \Bigr)\Bigl|\frac{\partial u_1}{\partial y_1}\Bigr|
		\]
		so that ${a(u_1,u_0)=-a(u_0,u_1)}$ and $c(u_1,u_0)=0=-c(u_0,u_1)$. Thus both cases satisfy the condition of Proposition \ref{propGroup}.
	\end{example}

	\begin{example}\label{exCircle2}
			For a rotation model on $\mathbb{R}^p$, $\mathscr{Y}\subset \mathbb{R}^p$, $\Gamma=[0,2\pi)$ and the group action is matrix multiplication by a rotation matrix $g\in G\subseteq \text{SO}(p)$, where $\text{SO}(p)$ is the special orthogonal group of dimension $p$. The identity element is matrix multiplication by the $p$-dimensional identity matrix and $f_U(u;\gamma)$ is a distribution on the hypersphere, where $\gamma$ is a rotation parameter. The von Mises-Fisher distribution is one such example. Consider the case of $p=2$.
	
	Let $v\in \mathbb{R}^2$ be a point on a circle of fixed arbitrary radius, either a location parameter of a probabilistic model or a data point on the circle. In a symmetric parametrization of a rotation model, $v_\gamma=g_\gamma v$ defines a (counter-clockwise) rotation of angle $v_\gamma^\T v = \gamma$ and $v_1 = g_\psi v_\gamma$, $v_0=g_\psi^{-1} v_\gamma$ define (counter-clockwise and clockwise) rotations of angles $v_1^\T v_\gamma=\psi$ and $v_0^\T v_\gamma=-\psi$ respectively, where
	\begin{eqnarray*}
		g_\psi^{-1}&=&\begin{pmatrix} \cos \psi & -\sin \psi \\ \sin \psi & \;\;\,\,\cos \psi \end{pmatrix}^{-1}=\begin{pmatrix} \;\;\,\,\cos \psi & \sin \psi \\ -\sin \psi & \cos \psi \end{pmatrix}\\
		&=&\begin{pmatrix} \cos (-\psi+2\pi) & -\sin (-\psi+2\pi) \\ \sin (-\psi+2\pi) & \;\;\,\,\cos (-\psi+2\pi) \end{pmatrix} \in G.
	\end{eqnarray*}
	and $g_\psi g_\gamma=g_{\gamma+\psi}\in G$, $g_\psi^{-1} g_\gamma=g_{\gamma-\psi}\in G$. Thus, in equation \eqref{eqSymmPairGeneral} $g\gamma$ should be understood either as addition of angles $g\gamma=\gamma+\psi$ or multiplication of rotation matrices $g\gamma = g_\psi g_\gamma$. 
	
	Let $y_1=(y_{11}, y_{12})^\T$, $y_0=(y_{01}, y_{02})^\T$. Then
	\[
	u_1 = g^{-1}y_1 = \begin{pmatrix}
	\;\;\,\,y_{11}\cos \psi + y_{12}\sin \psi \\
	-y_{11} \sin \psi +y_{12} \cos \psi
	\end{pmatrix}, \quad 	u_0 = g y_0 = \begin{pmatrix}
	y_{01}\cos \psi - y_{02}\sin \psi \\
	y_{01} \sin \psi +y_{02} \cos \psi
	\end{pmatrix}.
	\]
	The Jacobian determinants $J_{1}^+$ and $J_{0}^{+}$ are both equal to $|\cos^2\psi+\sin^2\psi|=1$. Thus, $c(u_1,u_0)=0=-c(u_0,u_1)$ in Proposition \ref{propGroup}. Since 
	\begin{eqnarray*}
		\frac{\partial u_1}{\partial \psi} \hspace{4pt} =& \hspace{-4pt} \begin{pmatrix}
			-y_{11}\sin \psi + y_{12}\cos \psi \\
			-y_{11} \cos \psi -y_{12} \sin \psi
		\end{pmatrix} \hspace{-4pt}& = \hspace{4pt} g_\theta^{-1} u_1, \\
		\frac{\partial u_0}{\partial \psi} \hspace{4pt} =& \hspace{-4pt} \begin{pmatrix}
			-y_{01}\sin \psi - y_{02}\cos \psi \\
			\;\;\,y_{01} \cos \psi -y_{02} \sin \psi
		\end{pmatrix} \hspace{-4pt}& = \hspace{4pt} g_\theta u_0,
	\end{eqnarray*}
	with $\theta=\pi/2$,	the quantity $a(u_1,u_0)$ from Proposition \ref{propGroup} is
	\[
	a(u_1,u_0)=\begin{pmatrix} \cos (-\pi/2) & -\sin (-\pi/2) \\ \sin (-\pi/2) & \;\;\,\,\cos (-\pi/2) \end{pmatrix} u_1 + \begin{pmatrix} \cos (\pi/2) & -\sin (\pi/2) \\ \sin (\pi/2) & \;\;\,\,\cos (\pi/2) \end{pmatrix} u_0 = -a(u_0,u_1),
	\]
	verifying the conditions for the $\psi$-symmetric parametrization of rotation families on the circle.
	\end{example}

	\section{Examples}\label{secExamples}

	\subsection{Matched pairs}\label{secMissRE}

	Let $Y_{i1}$ and $Y_{i0}$ be random variables corresponding to observations on treated and untreated individuals in a matched pair design, where $i=1,\ldots, n$ is the pair index. The treatment effect is represented by a parameter $\psi$ while the pair effects are encapsulated in the pair-specific nuisance parameter $\gamma_i$, avoiding explicit modelling assumptions in terms of covariates.  As noted in Example \ref{exRunningExample}, two approaches to analysis treat the pair effects as fixed arbitrary constants, or as independent and identically distributed random variables. In the latter case, it is common to model the distribution parametrically, typically producing efficiency gains in the treatment effect estimator if the parametric assumptions hold to an adequate order of approximation. 
	
	The relevant calculations are the same for every pair, and we therefore omit the pair index $i$ from the notation. Let the true joint density function of $(Y_{1},Y_{0})$ be given by
	\begin{equation}\label{eqEta}
	m(y_1,y_0)=\int f_1(y_1\hspace{0.5pt};\psi^*, \gamma)f_0(y_0\hspace{0.5pt};\psi^*, \gamma)f(\gamma)d\gamma
	\end{equation}
	with $f(\gamma)$ an unknown density function for the nuisance parameters, which are treated as independent and identically distributed across pairs.
	 
	Conditionally on $\gamma$, let $\ell(\psi;\gamma,y_1,y_0)$ denote the log-likelihood contribution for $\psi$ from a single pair. 
	
	\begin{proposition}\label{propConsistencyPairs}
		Suppose that the assumed model over $\gamma$ is parametrized by $\lambda$, producing a log-likelihood contribution $\ell(\psi,\lambda; y_1,y_0)$, assumed strictly concave. Suppose further that, conditionally on $\gamma$, $Y_1$ and $Y_0$ have a joint distribution that is parametrized $\psi$-symmetrically $($Definition \ref{defSymmetric}$)$. Then provided that the group induces antisymmetry on the log-likelihood derivative $($Definition \ref{defAntisymm}$)$, it follows that $\psi^* \perp_m \Lambda$ and
\begin{equation}\label{eqSuffLik}
0 = \mathbb{E}_m[\nabla_\psi \ell(\psi^*,\lambda)] = \int_{\mathscr{Y}}\int_{\mathscr{Y}} \bigl(\nabla_\psi \ell(\psi^*,\lambda \hspace{0.5pt}; y_1, y_0)\bigr)m(y_1,y_0)dy_1 dy_0
\end{equation}
		for all $\lambda\in\Lambda$. Thus, $\psi_m^0 = \psi^*$ by Proposition \ref{propConsistency}.
	\end{proposition}

	\begin{example}[Continuation of Examples \ref{exRunningExample} and \ref{exRunningEx2}]\label{exExpMP}
		Example \ref{exRunningEx2} establishes antisymmetry of the log-likelihood derivative for the symmetric parametrization of the exponential matched pairs problem. Consistency of $\hat\psi$ thus follows by Proposition \ref{propConsistencyPairs} regardless of the true distribution and assumed model for the random effects $\gamma$. This considerably extends the result of \cite{BC2020} and provides deeper insight into the structure of inference under model misspecification.
\end{example}
	
	\begin{remark}
It is not too difficult to check following calculations similar to that in Appendix A.4 of \cite{BC2020} that if we parametrize this model non-symmetrically, for example with rates $\gamma_i\theta$ and $\gamma_i$, the resulting estimator is not consistent if the gamma random effects assumption is erroneous.
\end{remark}

	Proposition \ref{propGroup} gives easily verifiable conditions for antisymmetry of the log-likelihood derivative, which holds for all location and scale families in Example \ref{exLocScale}. 
	The following explicit calculation for the normal distribution mirrors Example \ref{exRunningEx2}.
	
	\begin{example}[normal matched pairs with arbitrary mixing]\label{exNormalMP}
		
		Let $Y_1$ and $Y_0$ be normally distributed with unit variance and  means $\gamma+\psi^*$ and $\gamma-\psi^*$ respectively, $\gamma$ being treated as random. The conclusion is unchanged if an unknown variance parameter is present. Conditionally on $\gamma$, with $u_1 = y_1-\psi^*$, $u_0 = y_0 + \psi^*$, the joint density function of $(Y_1,Y_0)$ is
	\[
f_1(y_1 \hspace{0.5pt};\psi^*,\gamma)f_0(y_0 \hspace{0.5pt}; \psi^*,\gamma)dy_1dy_0=\exp[-\{(u_1-\gamma)^2 + (u_0-\gamma)^2\}/2]du_1du_0,
\]
		which is symmetric in $u_1$ and $u_0$ by construction as a result of using a $\psi$-symmetric parametrization with respect to addition, the relevant group action for location models. Conditionally on $\gamma$, the likelihood as a function of $\psi$ is 
	\begin{eqnarray}\label{eqNormalLik}
\nonumber 	L(\psi; y_0, y_1, \gamma)	& \propto & \exp\{-\tfrac{1}{2}(y_1-\gamma-\psi)^2 - \tfrac{1}{2}(y_0-\gamma+\psi)^2 \} \\ &\propto& \exp[-\tfrac{1}{2}\{(y_1-\psi)^2 +(y_0+\psi)^2\}],
\end{eqnarray}
		where the proportionality symbol on the second line absorbs the multiplicative term $\exp\{-\gamma^2+\gamma(y_1+y_0)\}$, which does not depend on $\psi$.
	    On taking logarithms, discarding additive terms independent of $\psi$ and differentiating, 
		\[
		\nabla_{\psi} \ell(\psi; u_0, u_1) =  y_1-y_0 -2\psi = u_1-u_0 
		= -\nabla_{\psi} \ell(\psi; u_1, u_0),
		\]
		verifying the antisymmetry condition of Definition \ref{defAntisymm}.
		
		For any assumed random effects distribution for $\gamma$ with density function parametrized by $\lambda$, an application of Proposition  \ref{propConsistencyPairs} guarantees that the resulting maximum likelihood estimator $\hat\psi$ is consistent for $\psi^*$ under arbitrary misspecification of the random effects distribution. 
	\end{example}
	
	\begin{remark}
		A separate point in connection with Example \ref{exNormalMP} is that, due to the separation in the log-likelihood function specific to this case (a parameter cut in the terminology of \cite{B-N1978}), $\hat{\psi}=(Y_1-Y_0)/2$ regardless of the assumed random effects distribution, and  for an analysis from $n$ pairs $\hat{\psi}=(S_1-S_0)/2n$ where $S_j=\sum_{i=1}^n Y_{ij}$ where $i$ is the pair index. For any finite-variance random effects distribution $f$, $\text{var}(Y_1)=\text{var}(Y_0)=1+\text{var}_f(\gamma)$ and $\text{cov}(Y_1,Y_0)=\text{var}_f(\gamma)$. Thus $\text{var}(\hat{\psi})=2/n$ irrespective of $f$ and of the assumed (erroneous) model over the random effects.
	\end{remark}

\begin{remark}
In the context of stratum-specific random effects, Lindsay \cite[Thm 4.2]{Lindsay1985} established the limiting distribution of his proposed estimator of the interest parameter, valid whether or not the random effects distribution is misspecified, under the assumption that his estimator is consistent. The latter aspect, central to the present paper, was not discussed in that work.	
\end{remark}

	\subsection{Two-group problems}\label{secTwoGroup}

	Section \ref{secMissRE} dealt with an experimental setting in which the matching ensures balance. A more realistic situation in an observational context is that observations on treated and untreated individuals are stratified into groups that are as similar as possible. This in general produces unbalanced strata, having a potentially different number of individuals in the treated and untreated groups for each stratum.
	
	Inference is based on the sufficient statistics $S_{j1}$ and $S_{j0}$ within treatment groups and strata. Let $(Y_{ij1})_{i=1}^{r_{j1}}$ and $(Y_{ij0})_{i=1}^{r_{j0}}$ be observations within the $j$th stratum for treated and untreated individuals respectively. Some examples fix ideas. If the observations $(Y_{ij1})_{i=1}^{r_{j1}}$ and $(Y_{ij0})_{i=1}^{r_{j0}}$ are normally distributed with means $\gamma_{j}+\psi^*$ and $\gamma_{j}-\psi^*$ and variance $\tau$, the likelihood contribution to the $j$th stratum depends on the data only through $\sum_{i=1}^{r_{j1}}Y_{ij1}/r_{j1}$ and $\sum_{i=1}^{r_{j0}}Y_{ij0}/r_{j0}$, whose distributions are normal with means unchanged, and variances $\tau/r_{j1}$ and $\tau/r_{j0}$ respectively. If, conditionally on $\gamma_{j}$, the individual observations are Poisson distributed counts of rates $\gamma_{j}\psi^*$ and $\gamma_{j}/\psi^*$, the sufficient statistics are sums of these counts, Poisson distributed of rates $r_{j1}\gamma_{j}\psi^*$ and $r_{j0}\gamma_{j}/\psi^*$. As a third example, if the originating variables are exponentially distributed of rates $\gamma_{j}\psi^*$ and $\gamma_{j}/\psi^*$, the sufficient statistics are gamma-distributed sums of shape and rate parameters $(r_{j1}, \gamma_j\psi^*)$ and $(r_{j0},\gamma_j/\psi^*)$ respectively. The matched comparison setting of \S \ref{secMissRE} is a special case of this formulation with $r_{j1}=r_{j0}=1$ for all $j$. The more general balanced case with $r_{j1}=r_{j0}= r$ is also implicitly covered by Propositions \ref{propGroup} and $\ref{propConsistencyPairs}$.
	
	The situation when $r_{j1}\neq r_{j0}$ is more complicated as it implies that even if the distributions of $S_{j1}$ and $S_{j0}$ belong to a group family, the transformations required to express their distributions in the standardized form $f_U$ depend on the different values of $r_{j1}$ and $r_{j0}$, so that Definition \ref{defSymmetric} is violated. The following extension of Proposition \ref{propConsistencyPairs}, which assumes the same random effects formulation for the nuisance parameters $\gamma_j$, omits the subscripts for the strata as before.
	
	\begin{proposition}\label{prophsymmetric}
		Suppose that, conditionally on $\gamma$, $S_1$ and $S_0$ are independent random variables with probability measures in $\mathcal{P}_G$. Let $h_1\in G$ and $h_0\in G$, not depending on $\psi$ or $\gamma$, be such that $T_1=h_1 S_1$ and $T_0=h_0 S_0$ have a joint distribution that is parametrized $\psi$-symmetrically in the sense of Definition \ref{defSymmetric} and such that the log-likelihood derivative, when expressed in terms of $u_1=g^{-1}t_1$ and $u_0=g t_0$, is antisymmetric in the sense of Definition \ref{defAntisymm}. Then, provided that $\ell(\psi,\lambda)$ is strictly concave, $\psi^*\perp \Lambda$ and 
		\begin{equation}\label{eqSuffLik2}
		0 = \mathbb{E}_m[\nabla_\psi \ell(\psi^*,\lambda)] = \int\int \bigl(\nabla_\psi \ell(\psi^*,\lambda \hspace{0.5pt}; s_1, s_0)\bigr)m(s_1,s_0)ds_1 ds_0
		\end{equation}
		for all $\lambda\in \Lambda$. Thus, $\psi_m^0 = \psi^*$ by Proposition \ref{propConsistency}.
	\end{proposition}
	
	For discrete problems such as the two-group Poisson situation mentioned above, there is no natural group and Propositions \ref{propConsistencyPairs} and \ref{prophsymmetric} are therefore inapplicable. Direct verification of Proposition \ref{propCoxWong} is sometimes more fruitful, as illustrated in the following examples.
	
	\begin{example}[stratified two-group Poisson problem with unbalanced strata]\label{exCW}
		
		Cox and Wong~\cite{CoxWong2010} considered misspecification of a random effects distribution in a stratified two-group problem. The distributions of counts $S_{j1}$ and $S_{j0}$ in the treated and untreated groups for stratum $j$ are, conditionally on $\gamma_j$, Poisson with means $r_{j1}\gamma_j \exp (\theta^*)$ and $r_{j0}\gamma_j \exp (-\theta^*)$ respectively, where $r_{j1}$ and $r_{j0}$ reflect the number of patients at risk in each group and $\gamma_j$ are stratum-specific nuisance parameters. 
		
		The problem can equivalently be parametrized in terms of $\psi=e^\theta$ but the above formulation is that used by Cox and Wong \cite{CoxWong2010}. The result of Proposition \ref{propCWPoisson} was noted there; we complete their proof in Appendix \ref{secProofs}.
		
		\begin{proposition}[Cox and Wong, \cite{CoxWong2010}]\label{propCWPoisson}
			Suppose that in Example \ref{exCW}, the nuisance parameters $(\gamma_j)_{j=1}^n$ are modelled as independent and identically distributed random variables with gamma density $\xi(\xi\gamma)^{\omega - 1}\exp(-\xi \gamma)/\Gamma(\omega)$, parametrized in terms of the notionally orthogonal parameters $\nu=\xi/\omega$ and $\omega>0$. Provided that there is a member of the assumed model that gives the same expectation $\EE_{\nu,\omega}(\gamma_j)=1/\nu$ as under the true random effects distribution for $(\gamma_j)_{j=1}^n$, $\hat{\theta}$ is consistent for $\theta^*$.
		\end{proposition}
	\end{example}

	\begin{example}[time to first event in a two-group problem with unbalanced strata]\label{exCWV2}
		
		We describe here the generating mechanism for an exponential analogue of Example \ref{exCW}. 
		Suppose the available data in Example \ref{exCW} consist also of the failure times. If these are, for each individual, independent exponentially distributed of rates $\gamma_j \exp(\theta^*)$ and  $\gamma_j \exp(-\theta^*)$ in the $j$th stratum, with $r_{j1}$ and $r_{j0}$ individuals at risk in the treated and untreated groups, then the times $S_{j1}$ and $S_{j0}$ to the first failure in each group are the group minima, exponentially distributed with rates $r_{j1}\gamma_j \exp(\theta^*)$ and $r_{j0}\gamma_j \exp(-\theta^*)$, respectively. By using only the minima across strata, information is sacrificed, which would normally not be advised unless to check on modelling assumptions. After reparametrization $\psi=e^\theta$, the problem satisfies the conditions of Proposition \ref{prophsymmetric} with $T_{j1}=r_{j1}S_{j1}$ and $T_{j0}=r_{j0}S_{j0}$. Thus $\hat{\psi}$ is consistent for $\psi^*=\exp(\theta^*)$.
	\end{example}

	\subsection{Generalized linear models}\label{secExpFamily}
	
	A simple and well known example in which Proposition \ref{propConsistency} is directly applicable is misspecified dispersion in exponential-family regression problems. Suppose that the true density function belongs to the class of canonical exponential family regression models with unknown regression coefficient $\psi$ on a set of known covariates $x_1,\ldots,x_n$ (each column vectors of the same dimension as $\psi$), and dispersion parameters $\phi_1,\ldots,\phi_n$. The form of the exponential-family log-likelihood function with respect to $\phi_i$ is such that $\psi^* \perp_{m} \phi_i$ for all $i$. In particular this implies that if $\phi_i$ is modelled erroneously as $\phi_i(\lambda)$, perhaps depending on covariates other than or including $x_i$, $\psi^* \perp_{m} \Lambda$. Write the log-likelihood function of the assumed model as 
	\begin{equation}\label{eqGLMlik2}
	\ell(\psi,\lambda) = \sum_{i=1}^{n}\phi_i^{-1}(\lambda)\{y_i x_i^\T\psi - K(x_i^\T\psi)\}+\sum_{i=1}^n h\{y_i; \phi_i(\lambda)\}.
	\end{equation}
	We do not distinguish notationally between realisations of the random variables $Y_1,\ldots,Y_n$ and arbitrary evaluation points for their joint density function conditional on covariates. With outcomes treated as random, the estimating equation for $\psi$ is unbiased in the sense that $\mathbb{E}_m[\nabla_\psi \ell(\psi^*,\lambda)]=0$ for all $\lambda\in\Lambda$. Thus $\hat{\psi}$ is consistent in the presence of arbitrary misspecification of the dispersion parameters. 
	
	In a continuation of this generalized linear models example, we assume this time a simplified specification for dispersion but allow for relevant covariates to be omitted. Let $Z=(X, W)$ denote the included and omitted covariates. The conditional probability density or mass function of $Y=(Y_{1},\ldots,Y_n)$ at $y=(y_1,\ldots,y_n)$ is, with $\theta^*=(\psi^{*\T},\lambda^{*\T})^\T$,
\[
m(y;z_{1}^\T\theta^*, \ldots, z_n^\T\theta^*)=\exp\Bigl[\phi^{-1}\Bigl\{\theta^{*\T}\sum_{i=1}^{n}z_i y_i - \sum_{i=1}^n K(z_i^{T}\theta^*)\Bigr\}\Bigr]\prod_{i=1}^n h(y_i, \phi),
\]
	Let $\check{m}(y;x_1^\T\psi,\ldots,x_n^\T\psi)$ be the assumed model, in which the outcome is modelled only in terms of $X$. The orthogonality condition $\psi^* \perp_m \Lambda$ is 
	\begin{equation}\label{eqOrthoGLM}
	\sum_{i=1}^n K^{\prime\prime}(x_i^\T \psi^* + w_i^\T \lambda)x_i w_i^\T =0
	\end{equation}
	for all $\lambda$, a highly restrictive condition. By Proposition \ref{propConsistency}, the additional conditions under which $\psi_m^0 = \psi^*$ are those that set $\EE_m[\nabla_\psi \ell(\psi^*,\lambda)]=0$ for all $\lambda$, where $\ell$ is the log-likelihood function under the assumed model. Since $\lambda$ is not present under the assumed model,
\[
\EE_m[\nabla_\psi \ell(\psi^*,\lambda)]
=  \sum_{i=1}^{n}\{K^\prime(z_{i}^\T\theta^*) - K^\prime(x_{i}^\T\psi^*)\} x_i = \sum_{i=1}^n K^{\prime\prime}(x_i^\T \psi^* + \overline{\xi}_i)x_i w_i^\T\lambda^*
\]
	where $\overline\xi_i$ is on a line between $0$ and $w_i^\T \lambda^*$ and the first equality is because $K^\prime(z_{i}^\T\theta^*) = \EE_{m}(Y_i\mid Z_i=z_i)$. Since $K^{\prime\prime}(x_i^\T \psi^* + \overline{\xi}_i)\neq K^{\prime\prime}(x_i^\T \psi^* + w_i^\T \lambda)$ in general, the second condition $\EE_m[\nabla_\psi \ell(\psi^*,\lambda)]=0$ typically does not hold exactly even when the columns of $X$ and $W$ are orthogonal. An exception is the normal-theory linear model, for which $K(\zeta)=\zeta^2/2$ so that $K^{\prime \prime}(\zeta)=1$. For an insightful discussion of non-collapsibility in logistic and other regression models, see \cite{Daniel}, and for some relevant approximations, see \cite{CoxWermuth1990}.
	
	Although both conclusions of this subsection are well established in the literature, their purpose is to illustrate application of Proposition \ref{propConsistency} to some well-known settings.
	
	\subsection{Marginal structural models}\label{secMSM}
	
	Evans and Didelez~\cite{ED2023} presented a further example in which consistency arises as a special case of Proposition \ref{propConsistency}. Their paper introduced the idea of a frugal parametrization of a marginal structural model, constructed to complete the specification of a model respecting the interventional components of interest, but without loss or redundancy. In a causal system of the form, e.g.,~$Z\rightarrow X$, $Z \rightarrow Y$, $X\rightarrow Y$, the authors suggest a reparameterization of the joint distribution $p_{ZXY}$ from $(p_Z, p_{X|Z}, p_{Y|ZX})$ to $(p_Z, p_{X|Z}, p^*_{Y|X}, \phi^*_{YZ|X})$ where $p^*_{Y|X}$, with finite-dimensional parameter $\theta^*_{Y|X}$, is typically the unobservable interventional distribution of interest for $Y \mid \text{do}(X)$, and $\phi^*_{YZ|X}$ is whatever is necessary to complete the joint model. The quantity $p_{X|Z}$ is  known as the propensity score. We refer to the original paper for a detailed discussion of these quantities, which in turn are expressed in terms of finite-dimensional parameters, say $\psi$ and $\lambda$. The vector parameter $\psi$ with true value $\psi^*$ comprises $\phi^*_{YZ|X}$ and the main component of interest $\theta^*_{Y|X}$ which together specify $p_Z p_{Y|ZX}$, while $\lambda$ characterizes the propensity score $p_{X|Z}$. A consequence of the above parametrization within the class of marginal structural models is Theorem 5.1 of \cite{ED2023}, which establishes consistency and asymptotic normality of the maximum likelihood estimator $\hat{\psi}$ of $\psi^*$, even if $p_{X|Z}=p_{X|Z}(\lambda)$ is misspecified. As noted in \cite{ED2023}, the conclusion arises as a consequence of the parameter cut between $p_Z p_{Y|ZX}$ and $p_{X|Z}$, and therefore between $\psi$ and $\lambda$. The parameter cut  implies the parameter orthogonality condition $\psi^*\perp_m \Lambda$, and \cite{ED2023} implicitly establish in their proof  the remaining condition of Proposition \ref{propConsistency}.

	\section{Discussion}\label{secDiscussion}
	
	\subsection{Stronger inferential guarantees}\label{secInf}
	
	Even if consistency holds, inferential guarantees entail estimation of variance, a difficult problem due to failure of Bartlett's second identity. Such difficulties persist under the stronger condition of Proposition \ref{propOrthog}.

	In the notation of \S \ref{secConsistencyGeneral}, let 
	\[
	q = \mathbb{E}_{m}[(\nabla_{(\psi,\lambda)} \ell)( \nabla_{(\psi,\lambda)} \ell)^{\T}], \quad 
	\check{q}= \mathbb{E}_{(\psi,\lambda)}[(\nabla_{(\psi,\lambda)} \ell)( \nabla_{(\psi,\lambda)} \ell)^{\T}],
	\]
	both functions of $\psi$ and $\lambda$. Proposition \ref{propOrthog}, equation \eqref{eqLimit} and consistency of $\hat{\psi}$ imply
	\[
	\mathbb{E}_{(\psi,\lambda)}[\nabla_{(\psi,\lambda)}\ell(\psi^*, \lambda_{m}^0)] = 0,
	\]
	and differentiation gives Bartlett's second identity under the assumed model:
	\[
	-\check{\imath}(\psi^*, \lambda_{m}^0) + \check{q}(\psi^*, \lambda_{m}^0) = 0.
	\]
	The relevant quantities for establishing the limiting distributions of standard test statistics are, however, $i$ and $q$, where $i$ was defined in \S \ref{secConsistencyGeneral}. These do not in general satisfy $-i(\psi^*, \lambda_{m}^0)+q(\psi^*, \lambda_{m}^0)=0$, even with the additional structure of Proposition \ref{propOrthog}. An exception is Example \ref{exNormalMP} because of the particular simplicity of this case.
	
	It can be shown using now standard arguments that, with the log-likelihood function constructed from $n$ independent observations, provided that $\hat{\psi}$ is consistent for $\psi^*$,
	\begin{equation}\label{eqSandwich}
	\sqrt{n}(\hat{\psi}-\psi^*) \rightarrow_d N(0, (i^{-1}\hspace{0.5pt} q \hspace{0.5pt} i^{-1})_{\psi\psi}),
	\end{equation}
	the variance being the so-called sandwich formula, in which the constituent terms are evaluated at $(\psi^*, \lambda_{m}^0)$. This reduces to $i^{\psi\psi}(\psi^*, \lambda_{m}^0)=i^{\psi\psi}(\psi^*, \lambda^*)$ when the model is correctly specified, so that the usual asymptotic result is recovered. The standard references \cite{Cox1961}, \cite{Cox1962}, \cite{Kent1982}, \cite{White1982a}, and \cite{White1982b} have $\psi_m^0$ in place of $\psi^*$ in \eqref{eqSandwich}. Qualitatively similar results hold for the profile score and profile likelihood ratio statistics, the main point being that confidence set estimation is infeasible without knowledge of $m$ unless $\check{q}(\psi^*, \lambda_{m}^0)=q(\psi^*, \lambda_{m}^0)$ to some adequate order of approximation. In general this requires, not only that the expectations of the sufficient statistics are robust to misspecification, but also that the expectations of any relevant squares and cross terms are stable.
	
	Suppose, however, that the parametrization is orthogonal in the sense $\psi^* \perp_m \lambda_m^0$. Then
	\[
	(i^{-1}\hspace{0.5pt} q \hspace{0.5pt} i^{-1})_{\psi\psi} =  i^{\psi\psi} q_{\psi\psi} i^{\psi\psi},
	\]
	and the familiar result is also recovered if $q_{\psi\psi}i^{\psi\psi}=\text{I}_{\text{dim}(\psi)}$, a weaker requirement than $q=i$ especially if $\psi$ is a scalar. Under the additional structure of Proposition \ref{propOrthog}, inference based on the assumed model is unaffected by the misspecification.

	\subsection{Treating incidental parameters as fixed or random}\label{secFixedRandom}
	
	In the context of \S \ref{secMissRE}, the lack of general inferential guarantees beyond consistency may well be an argument for treating the nuisance parameters as fixed. Although the conceptual distinction is consequential for the analysis, a formulation in which incidental parameters are treated as fixed and arbitrary, and one in which they are treated as independent and identically distributed from a totally unspecified distribution, are numerically indistinguishable. In this sense, the fixed-parameter formulation is essentially nonparametric for the nuisance component, with the distinguishing feature that it evades estimation of the infinite-dimensional nuisance parameter.
	
	The group structure used to define the symmetric parametrization in Proposition \ref{propConsistencyPairs} also allows elimination of nuisance parameters by suitable preliminary manoeuvres when the nuisance parameters are treated as fixed. In the exponential matched pairs setting of Example \ref{exExpMP}, \cite{BC2020} compared  inference based on the distribution of $Y_{i1}/Y_{i0}$ to that based on modelling the pair effects by a gamma distribution. The random pair effects model is more efficient when the gamma distribution is correct, but efficiency degrades substantially under misspecification, compared to maximum likelihood estimation based on the distribution of ratios. 
	
	Proposition \ref{propConsistencyPairs} hinges on a symmetric parametrization having been chosen, which is possible only under the group structure of \S \ref{secSymmetry}. Outside of this setting there are no guarantees of consistency of the maximum likelihood estimator with misspecified random effects distribution, while it may still be possible to eliminate the pair effects with exact or approximate conditioning arguments, as illustrated in the second example of \cite{BC2020}.

	\subsection{Overstratification and other encompassing models}\label{secOverstratification}
	
	Two essentially equivalent ways to adjust for potential confounders are to stratify by observed explanatory features, leading to strata effects $(\gamma_j)_{j=1}^m$ as in Example \ref{exCW}, or to adjust for the explanatory features in a regression analysis. Both approaches are a form of conditioning by model formulation. It is arguably more relevant in this context to treat the strata effects as fixed rather than random. In De Stavola and Cox \cite{deStavola} this modification of Example \ref{exCW} was considered, showing the extent to which overstratification decreases efficiency of the estimator. As is intuitively clear, no bias is incurred through overstratification, yet the conditions of Propositions \ref{propConsistency} and \ref{propCoxWong} do not always hold, as illustrated by the following example. The explanation is that the true density is nested within the encompassing model, so that $\lambda_m^0=\lambda^*=0$. The overstratified model is therefore not misspecified according to the definition below \eqref{eqLimit}.
	
	To isolate the point at issue, we assume that any  dispersion parameters are equal to 1, the conclusion of Proposition \ref{propOverstrat} being unchanged for arbitrary dispersion parameter by the discussion of \S \ref{secExpFamily}.
	
	\begin{proposition}\label{propOverstrat}
		Suppose that, conditional on observed explanatory features, the probability density or mass function of $Y=(Y_{1},\ldots,Y_n)$ at $y=(y_1,\ldots,y_n)$ is of canonical exponential family regression form, that is,
	\[
m(y;x_{1}^\T\psi^*, \ldots, x_n^\T\psi^*)=\exp\Bigl\{\sum_{i=1}^{n} y_i x_i^\T\psi^{*}  - \sum_{i=1}^n K(x_i^{T}\psi^*)\Bigr\}\prod_{i=1}^n h(y_i).
\]	
		Suppose further that the analysis is overstratified so that under the assumed model the log-likelihood function is $\ell(\psi,\lambda) = s^\T\psi + t^\T\lambda - K(\psi,\lambda)$, where $s=\sum_{i=1}^n x_i y_i$, $t=\sum_{i=1}^n w_i y_i$ and $K(\psi,\lambda)=\sum_{i=1}^n K(x_i^\T\psi + w_i^\T\lambda)$. The conditions of Propositions \ref{propConsistency} and \ref{propCoxWong} are in general violated, yet the maximum likelihood estimator $\hat{\psi}$ is consistent for $\psi^*$.
	\end{proposition}
	
	Another type of encompassing model arises in random effects models when the assumed family of distributions for the random effects is rich enough to include the true distribution, for instance if the postulated mixture class is nonparametric (Kiefer \& Wolfowitz, \cite{KW1956}), as discussed in \S \ref{secFixedRandom}. This type of situation can be understood in terms of the convex geometry of mixture distributions (Lindsay, \cite{Lindsay1983}) and was specialized to the case of binary matched pairs with logistic probabilities by Neuhaus et al.~\cite{Neuhaus1994}. Their conditions on the mixing distributions with regards to the resulting cell probabilities effectively imply that the model is correctly specified according to the definition below \eqref{eqLimit}. 
	
	\subsection{Estimating equations and Neyman orthogonality}\label{secNeyman}
	
	A generalization of the score equation \eqref{eqLimit} is to other estimating equations, often but not always obtained as the vector of partial derivatives of a convex loss function. The model need only be partially specified in terms of a known function $h$ of the data, an interest parameter $\psi$ and a nuisance parameter $\lambda$, such that the expectation when the assumed model is true satisfies $\mathbb{E}_{m}[h(\psi^*,\lambda^*;Y)]=\mathbb{E}_{(\psi^*,\lambda^*)}[h(\psi^*,\lambda^*;Y)]=0$. If the model is misspecified the limiting solutions satisfy, in analogy with \eqref{eqLimit},
	\[
	\mathbb{E}_{m}[h(\psi_m^0,\lambda_m^0; Y)]=0.
	\]
	A special choice of $h$, establishing a connection to the double robustness literature, is the Neyman-orthogonal score for $\psi$, defined as
	\[
	s_{\text{N}}^*(\psi;\lambda,Y) = \ell_\psi(\psi,\lambda) - w^{*\T} \ell_\lambda(\psi,\lambda),
	\]
	where $\ell_\psi$ and $\ell_\lambda$ are the partial derivates of $\ell$ with respect to $\psi$ and $\lambda$ and $w^{*\T}$ is a matrix of dimension $\text{dim}(\psi)\times \text{dim}(\lambda)$ given by $w^{*\T}=i^{*}_{\psi\lambda}i^{*-1}_{\lambda\lambda}$ when the model is correctly specified. Here $i^*_{\psi\lambda}$ and $i^*_{\lambda\lambda}$ are the components of the Fisher information matrix at the true parameter values $(\psi^*,\lambda^*)$. Consider instead 
	\begin{equation}\label{eqNOScore}
	s_{\text{N}}(\psi;\lambda,Y) = \ell_\psi(\psi,\lambda) - w^{\T} \ell_\lambda(\psi,\lambda), \quad w^{\T}=i_{\psi\lambda}i^{-1}_{\lambda\lambda},
	\end{equation}
	where $i_{\psi\lambda}$ and $i_{\lambda\lambda}$ are as defined in \eqref{eqInv}. Write the condition $i^{\psi\psi}g_\psi + i^{\psi\lambda}g_\lambda=0$ of Proposition \ref{propCoxWong} as
	\begin{equation}\label{eqCondNO1}
	g_\psi + i_{\psi\psi.\lambda}i^{\psi\lambda}g_\lambda=0, \quad  i_{\psi\psi.\lambda}=(i^{\psi\psi})^{-1}=i_{\psi\psi}-i_{\psi\lambda}i_{\lambda\lambda}^{-1}i_{\lambda\psi}
	\end{equation}
	where $i_{\psi\psi.\lambda}$ is the Fisher information for $\psi$ at $(\psi^*,\lambda)$, computed under the true model, having adjusted for estimation of $\lambda$. On noting that $i^{\psi\lambda}=-i^{\psi\psi}i_{\psi\lambda}i_{\lambda\lambda}^{-1}$, we see that the information identity $i^{\psi\psi}g_\psi + i^{\psi\lambda}g_\lambda=0$ of Proposition \ref{propCoxWong} is equivalent to requiring that the Neyman orthogonal score \eqref{eqNOScore}, with score adjustment $w$ computed under the true distribution, has zero expectation under the same distribution for all $\lambda\in \Lambda$. Specification of $w^\T$ in \eqref{eqNOScore} when the model is misspecified thus hinges on the conditions of Proposition \ref{propOrthog} being satisfied.

	\begin{appendix}
		\section{Detailed derivations}\label{secProofs}
		
		\subsection{Proof of Proposition \ref{propConsistency}}
		
		\begin{proof}
			By definition $\psi_m^0$ solves $\mathbb{E}_{m}[\nabla_{(\psi,\lambda)}\ell(\psi_m^0, \lambda_m^0)]=0$ and therefore also $\mathbb{E}_{m}[\nabla_{\psi}\ell(\psi_m^0, \lambda_m^0)]=0$. Taylor expansion of $\nabla_{\psi}\ell(\psi,\lambda_m^0)$ around $\psi^*$ for fixed $\lambda_m^0$ followed by evaluation at $\psi_m^0$ gives
			\[
			\psi_m^0 - \psi^* = - \overline{B}^{-1}\mathbb{E}_m[\nabla_{\psi}\ell(\psi^*,\lambda_m^0)],
			\]
			where, by the mean value theorem for vector valued functions,
			\[
			\overline{B}=\mathbb{E}_m \Bigl[\int_{0}^1 \nabla_{\psi\psi}^2 \ell (t\psi_m^0+(1-t)\psi^*, \lambda_{m}^0)dt\Bigr],
			\]
			with the integral taken elementwise. Since $\overline{B}$ is positive definite, $\psi_m^0 =\psi^*$ if and only if $\mathbb{E}_m\bigl[\nabla_{\psi}\ell(\psi^*,\lambda_m^0)\bigr]=0$. 
			
			A Taylor expansion of $\nabla_{\psi}\ell(\psi^*,\lambda)$ around $\lambda_m^0$ for fixed $\psi^*$ shows that
			\begin{equation}\label{eqExpAlpha}
			\mathbb{E}_m[\nabla_{\psi}\ell(\psi^*,\lambda)] - \mathbb{E}_m[\nabla_{\psi}\ell(\psi^*,\lambda_m^0)] = \overline{A} (\lambda - \lambda_m^0),
			\end{equation}
			where 
			\[
			\overline{A}=\mathbb{E}_m \Bigl[\int_{0}^1 \nabla_{\psi\lambda}^2 \ell (\psi^*, t\lambda+ (1-t)\lambda_m^0)dt\Bigr]=\int_{0}^1\mathbb{E}_m \bigl[ \nabla_{\psi\lambda}^2 \ell (\psi^*, t\lambda+(1-t)\lambda_m^0)\bigr]dt.
			\]
			Nullity of the left hand side of \eqref{eqExpAlpha} for general $\lambda$ entails $\overline{A}=0$, for which a sufficient condition is $\psi^* \perp_m \Lambda$. Without such orthogonality, positive and negative parts must cancel in the integral defining $\overline{A}$, which cannot be simultaneously achieved for all $\lambda$, establishing necessity of $\psi^* \perp_m \Lambda$ if $
			\mathbb{E}_{m}\bigl[\nabla_{\psi}\ell(\psi^*, \lambda)\bigr]=0$ is to hold for all $\lambda$.
		\end{proof}
		
		\subsection{Proof of Proposition \ref{propCoxWong}}
		
		\begin{proof}
			It is notationally convenient to suppose that the interest parameter $\psi$ is scalar, although the results generalize straightforwardly to higher dimensions. We therefore write $\ell_\psi$ for the derivative of the log likelihood function with respect to $\psi$ and $\ell_s$ for the derivative with respect to the $s$th component of $\lambda$, and similarly for other quantities. We use the convention that Roman letters appearing both as subscripts and superscripts in the same product are summed. Thus, for the purpose of this calculation, write $(\hat{\lambda}-\lambda)^s$ for the $s$th component of $\hat{\lambda}-\lambda$. 
			
			Let $\overline{\ell}=\ell/n$ and $\hat{\ell}=\overline{\ell}(\hat{\psi},\hat{\lambda})$, etc. The rescaling by $n$ is immaterial and is used only to aid presentation of the derivation. Taylor expansion of $\overline{\ell}_\psi$ around $(\psi^*, \lambda)$ followed by evaluation at $(\hat{\psi},\hat{\lambda})$ gives
			\begin{eqnarray*}
				0 = \hat{\ell}_\psi &=& \overline{\ell}_\psi + \overline{\ell}_{\psi\psi}(\hat{\psi}-\psi^*) + \overline{\ell}_{\psi s}(\hat{\lambda}-\lambda)^{s} \\
				& &  + \quad  \overline{\ell}_{\psi\psi t}(\hat{\psi}-\psi^*)(\hat{\lambda}-\lambda)^t 
				+ \tfrac{1}{2} \overline{\ell}_{\psi s t}(\hat{\lambda}-\lambda)^s (\hat{\lambda}-\lambda)^t + \cdots .
			\end{eqnarray*}
			
			Let $n \overline{j}^{rs}$ denote the components of the inverse matrix of the observed information matrix at $(\psi^*, \lambda)$ whose components are $n \overline{j}_{rs}=-\ell_{rs}$. Inversion of the previous equation as in Barndorff-Nielsen and Cox (1994, p.~149) gives
			\begin{eqnarray}\label{eqExpansion}
			\hat{\psi} - \psi^* &=&  \overline{j}^{\psi\psi}\overline{\ell}_\psi + \overline{j}^{\psi s}\overline{\ell}_s + \overline{j}^{\psi\psi}\overline{\ell}_{\psi\psi u}(\hat{\psi}-\psi^*)(\hat{\lambda}-\lambda)^u \\ 
			\nonumber & &  + \quad  \tfrac{1}{2}\overline{j}^{\psi \psi}\overline{\ell}_{\psi t u}(\hat{\lambda}-\lambda)^t (\hat{\lambda}-\lambda)^u + \tfrac{1}{2}\overline{j}^{\psi s}\overline{\ell}_{s t u}(\hat{\lambda}-\lambda)^t (\hat{\lambda}-\lambda)^u + \cdots.
			\end{eqnarray}
			Write $\overline{j}^{rs}=\overline{i}^{rs}+Z_{rs}/\sqrt{n}$, where $\overline{i}^{rs}$ denotes the components of the inverse of the Fisher information matrix $i_{rs}$ at $(\psi^*,\lambda)$ rescaled by $n$, and $Z_{rs}$ are  random variables of zero mean and $O_{p}(1)$ as $n\rightarrow \infty$. Similarly write $\overline{\ell}_s = \overline{g}_s + V_s/\sqrt{n}$ where $g_s$ is the expectation of $\overline{\ell}_s$. Expectations in $\overline{i}_{rs}$ and $\overline{g}_s$ are taken under the true model. It follows that the leading order term in the expansion is
			\[
			\hat{\psi} - \psi^* = \overline{i}^{\psi\psi}\overline{g}_\psi + \overline{i}^{\psi s}\overline{g}_s + \overline{i}^{\psi\psi} \frac{V_\psi}{\sqrt{n}} + \overline{i}^{\psi s} \frac{V_s}{\sqrt{n}} + \frac{Z_{\psi\psi}}{\sqrt{n}} \overline{g}_\psi + \frac{Z_{\psi s}}{\sqrt{n}} \overline{g}_s + \text{Rem}.
			\]
			If $\lambda=\lambda_m^0$, $\overline{i}^{\psi\psi}\overline{g}_\psi + \overline{i}^{\psi s}\overline{g}_s = O(n^{-1/2})$ provided that the dimension of $\lambda$ is treated as fixed. This is not in general the case when $\lambda\neq \lambda_m^0$, meaning that successive replacements in equation \eqref{eqExpansion} do not produce terms of decreasing orders of magnitude as in equation (5.19) of Barndorff-Nielsen and Cox (1994). 
			
			Suppose that 
			\begin{equation}\label{eqConditionSymm}
			\overline{i}^{\psi\psi}\overline{g}_\psi + \overline{i}^{\psi s}\overline{g}_s = 0
			\end{equation} 
			for any $\lambda$. Again, the second and higher order terms in \eqref{eqExpansion} do not converge to zero in general. The restriction to $\lambda = \lambda_{m}^0 + O(n^{-1/2})$ ensures that $\hat{\psi}$ converges to $\psi_m^0=\psi^*$ at the same rate as if the model was correctly specified. But if \eqref{eqConditionSymm} holds for all $\lambda$, it holds a fortiori for $\lambda$ in a neighbourhood of $\lambda_m^0$, thereby proving the claim.

			In the special case that $\psi$ and $\lambda$ are both scalar parameters, $i^{\psi\psi}=\text{det}^{-1}i_{\lambda\lambda}$ and $i^{\psi\lambda}=-\text{det}^{-1}i_{\psi\lambda}$, where $\text{det}=i_{\psi\psi}i_{\lambda\lambda}-i_{\psi\lambda}^2$, so that the general condition becomes $i_{\lambda\lambda} g_\psi = i_{\psi\lambda} g_\lambda$.

		\end{proof}

		\subsection{Proof of Proposition \ref{propGroup}}
		
		\begin{proof} Without loss of generality, write $\nabla_u f_U(u;\gamma)=b(u;\gamma)f_U(u;\gamma)$, which may entail defining $b(u;\gamma)$ trivially as $b(u;\gamma)=\nabla_u f_U(u;\gamma)/f_U(u;\gamma)$. The form of $b(u;\gamma)$ is immaterial for the argument to be presented. 
			
			The log-likelihood derivative is antisymmetric if and only if the likelihood derivative is antisymmetric, as $\nabla_\psi \ell = \nabla_\psi L/L$ and $L$ is symmetric in $u_1$ and $u_0$ by Definition \ref{defSymmetric}. Hence consider
			\begin{equation}\label{eqLikDeriv}
			\nabla_\psi L(\psi;\gamma,y_1,y_0)=  f_0(y_0;g^{-1}\gamma) \nabla_\psi f_1(y_1;g\gamma) + f_1(y_1;g\gamma)\nabla_\psi f_0(y_0;g^{-1}\gamma)
			\end{equation}
			where, on writing $u_1=g^{-1}y_1$ and $u_0=g y_0$,
			\[
			\nabla_\psi f_1(y_1;g\gamma) =    \Bigl|\frac{\partial u_1}{\partial y_1}\Bigr|\frac{\partial}{\partial \psi} f_U(g^{-1}y_1; \gamma) + f_U(g^{-1}y_1; \gamma) \frac{\partial}{\partial \psi}\Bigl|\frac{\partial u_1}{\partial y_1}\Bigr|
			\]
			and
			\[
			\frac{\partial}{\partial \psi} f_U(g^{-1}y_1; \gamma) = \frac{\partial}{\partial u_1} f_U(u_1; \gamma)\frac{\partial u_1}{\partial \psi} = b(u_1;\gamma)f_U(u_1;\gamma)\frac{\partial u_1}{\partial \psi}.
			\]
			Analogous expressions hold for $\nabla_\psi f_0(y_0;g^{-1}\gamma)$, giving
			\[
			\nabla_\psi f_1(y_1;g\gamma)=f_U(u_1;\gamma)\Bigl(b(u_1;\gamma)\frac{\partial u_1}{\partial\psi}\Bigl|\frac{\partial u_1}{\partial y_1}\Bigr| + \frac{\partial}{\partial \psi} \Bigl|\frac{\partial u_1}{\partial y_1}\Bigr| \Bigr)
			\]
			and
			\[
			\nabla_\psi f_0(y_0;g^{-1}\gamma)=f_U(u_0;\gamma)\Bigl(b(u_0;\gamma)\frac{\partial u_0}{\partial \psi}\Bigl|\frac{\partial u_0}{\partial y_0}\Bigr| + \frac{\partial}{\partial\psi} \Bigl|\frac{\partial u_0}{\partial y_0}\Bigr| \Bigr)
			\]
			respectively. On substituting in \eqref{eqLikDeriv} and expressing the resulting quantities in terms of $u_1$ and $u_0$ we obtain
			\[
			\nabla_\psi L(\psi;\gamma,u_1,u_0) = f_U(u_1;\gamma)f_U(u_0;\gamma)B(u_1,u_0)
			\]
			with
			\[
			B(u_1,u_0)=\Bigl(b(u_1;\gamma)\frac{\partial u_1}{\partial\psi}\Bigl|\frac{\partial u_1}{\partial y_1}\Bigr| + \frac{\partial}{\partial \psi} \Bigl|\frac{\partial u_1}{\partial y_1}\Bigr| \Bigr)\Bigl|\frac{\partial u_0}{\partial y_0}\Bigr| +	\Bigl(b(u_0;\gamma)\frac{\partial u_0}{\partial \psi}\Bigl|\frac{\partial u_0}{\partial y_0}\Bigr| + \frac{\partial}{\partial\psi} \Bigl|\frac{\partial u_0}{\partial y_0}\Bigr| \Bigr)\Bigl|\frac{\partial u_1}{\partial y_1}\Bigr|.
			\]
			so that the parametrization induces antisymmetry on the likelihood derivative in the sense of Definition \ref{defAntisymm} if and only if $B(u_1,u_0)=-B(u_0,u_1)$, which is seen to be equivalent to the conditions in Proposition \ref{propGroup} on noting that $J_1^+ J_0^+ = 1$ by the symmetry of the parametrization. 
		\end{proof}
		
		\subsection{Proof of Proposition \ref{propConsistencyPairs}}
		
		\begin{proof}
			The following calculations derive the conditions under which the sufficient and necessary condition from Proposition \ref{propConsistency} is satisfied, namely $\psi^* \perp_m \Lambda$ and
			\begin{equation}\label{eqSuffLik}
			0 = \mathbb{E}_m[\nabla_\psi \ell(\psi^*,\lambda)] = \int_{\mathscr{Y}}\int_{\mathscr{Y}} \bigl(\nabla_\psi \ell(\psi^*,\lambda \hspace{0.5pt}; y_1, y_0)\bigr)m(y_1,y_0)dy_1 dy_0,
			\end{equation}
			where 
			\begin{equation}\label{eqEta}
			m(y_1,y_0)=\int f_1(y_1\hspace{0.5pt};\psi^*, \gamma)f_0(y_0\hspace{0.5pt};\psi^*, \gamma)f(\gamma)d\gamma
			\end{equation}
			with $f(\gamma)$ an unknown density function for the nuisance parameter. All quantities are, for present purposes, evaluated at the same value of $\psi$ and the superscript is henceforth suppressed.
			
			Equation \eqref{eqSuffLik} is established first. Since $ \nabla_\psi \ell(\psi,\lambda)$ does not depend on $\gamma$ the order of integration in \eqref{eqSuffLik} can be interchanged, giving the sufficient condition
			\begin{equation}\label{eqSuffMP}
			\int_{\mathscr{Y}}\int_{\mathscr{Y}} \bigl(\nabla_\psi\ell(\psi,\lambda\hspace{0.5pt}; y_1, y_0)\bigr)f_1(y_1\hspace{0.5pt};\psi, \gamma)f_0(y_0\hspace{0.5pt};\psi, \gamma)dy_1 dy_0 = 0
			\end{equation}
			identically in $\gamma$. Suppose that, conditionally on $\gamma$, $Y_1$ and $Y_0$ have a joint distribution that is parametrized symmetrically in the sense of Definition \ref{defSymmetric}. On the left hand side of \eqref{eqSuffMP}, change variables to $u_1=g^{-1}y_1$ and $u_0 = g y_0$. The symmetry of the parametrization gives $J_1^+ J_0^+ = 1$ in Definition \ref{defSymmetric} and similarly $(\partial u_1/\partial y_1)^{-1}(\partial u_0/\partial y_0)^{-1}=1$. Therefore the volume element satisfies $dy_1 dy_0 = du_1 du_0$ and equation \eqref{eqSuffMP} is
			\begin{equation}\label{eqSuffMP2}
			\int_{\mathscr{Y}}\int_{\mathscr{Y}} \nabla_\psi\ell(\psi,\lambda \hspace{0.5pt}; u_1, u_0)f_U(u_1 \hspace{0.5pt}; \gamma)f_U(u_0 \hspace{0.5pt};\gamma)du_1 du_0 = 0,
			\end{equation}
			where the limits of integration are unaltered as a result of the assumed group structure. Since $f_U(u_1; \gamma)f_U(u_0;\gamma)du_1 du_0$ is symmetric in $u_1$ and $u_0$ and since $f_U$ is non-negative, a necessary and sufficient condition for equation \eqref{eqSuffMP2} is antisymmetry of $\nabla_\psi\ell(\psi,\lambda; u_1, u_0)$, in the sense
			\begin{equation}\label{eqAntiSymm}
			\nabla_\psi\ell(\psi,\lambda \hspace{0.5pt}; u_1, u_0)=-\nabla_\psi\ell(\psi,\lambda \hspace{0.5pt}; u_0, u_1).
			\end{equation}
			Here the derivative with respect to $\psi$ is taken first, followed by evaluation at $u_1=g^{-1}y_1$ and $u_0=gy_0$.
			
			Next consider the orthogonality condition of Proposition \ref{propConsistency}. The global orthogonality equation in the matched comparison problem is
			\[
			\int_{\mathscr{Y}}\int_{\mathscr{Y}} \nabla_{\psi\lambda}^2\ell(\psi,\lambda \hspace{0.5pt};y_1,y_0) m(y_1,y_0) dy_1 dy_0 = 0,
			\]
			where $m(y_1,y_0)=m(y_1,y_0 \hspace{0.5pt};\psi^*)$ is as defined in \eqref{eqEta} and $\psi^*$ has been reintroduced for clarity. Since 	$m$ does not depend on $\lambda$,
			\[
			\int_{\mathscr{Y}}\int_{\mathscr{Y}} \nabla_{\psi\lambda}^2\ell(\psi,\lambda \hspace{0.5pt};y_1,y_0) m(y_1,y_0) dy_1 dy_0 = \nabla_{\lambda} \int_{\mathscr{Y}}\int_{\mathscr{Y}} \nabla_\psi\ell(\psi,\lambda \hspace{0.5pt};y_1,y_0) m(y_1,y_0 \hspace{0.5pt};\psi^*) dy_1 dy_0.
			\]
			The double integral on the right hand side is, apart from the evaluation at $\psi^*$, that of \eqref{eqSuffLik}. Orthgonality at $\psi^*$ thus follows by an identical argument to that establishing \eqref{eqSuffLik}. The orthogonality is therefore local at $\psi^*$ and global in $\lambda$, as required in Proposition \ref{propConsistency}.
			
			The log-likelihood contribution from a single pair, treated as a function of both $(\psi,\lambda)$ and $(y_1,y_0)$ is
			\[
			\ell(\psi, \lambda \hspace{0.5pt}; y_1, y_0) = \log \int f_1(y_1 \hspace{0.5pt};\psi, \gamma)f_0(y_0 \hspace{0.5pt};\psi, \gamma)h(\gamma \hspace{0.5pt}; \lambda)d\gamma = \log L(\psi, \lambda \hspace{0.5pt}; y_1, y_0),
			\]
			say. The same argument as above shows that the likelihood function $L(\psi, \lambda \hspace{0.5pt}; y_1, y_0)$ is symmetric when expressed in terms of $u_1=g^{-1}y_1$ and $u_0=gy_0$. Since $\nabla_\psi\ell = \nabla_\psi L/L$, the antisymmetry condition \eqref{eqAntiSymm} is satisfied if $\nabla_\psi L$ is antisymmetric. The latter is implied by the clause of Proposition \ref{propConsistencyPairs} as
			\[
			\nabla_\psi L (\psi,\lambda \hspace{0.5pt}; y_1,y_0) = \int 	\nabla_\psi L (\psi\hspace{0.5pt}; \gamma, y_1,y_0)h(\gamma \hspace{0.5pt} ; \lambda)d\gamma.
			\]
		\end{proof}
		
		\subsection{Proof of Proposition \ref{prophsymmetric}}
		
		\begin{proof}
			
			The analogue of equation \eqref{eqSuffMP} in the proof of Proposition \ref{propConsistencyPairs} is
			\begin{equation}\label{eqSuffUnbalanced}
			\int_{\mathscr{Y}}\int_{\mathscr{Y}} \nabla_\psi\ell(\psi,\lambda\hspace{0.5pt}; s_1, s_0)f_1(s_1\hspace{0.5pt};\psi, \gamma)f_0(s_0\hspace{0.5pt};\psi, \gamma)ds_1 ds_0 = 0
			\end{equation}
			identically in $\gamma$. 
			
			In the likelihood derivative \eqref{eqLikDeriv}, the relevant terms are
			\begin{eqnarray*}
				\nabla_\psi f_1 (s_1;h_1 g \gamma) &=& \Bigl|\frac{\partial t_1}{\partial s_1}\Bigr|\nabla_\psi f_{T_1}(t_1;g\gamma) \\
				\nabla_\psi f_0 (s_0;h_0 g^{-1} \gamma) &=& \Bigl|\frac{\partial t_0}{\partial s_0}\Bigr|\nabla_\psi f_{T_0}(t_0;g^{-1}\gamma),
			\end{eqnarray*}
			and similarly for $f_1 (s_1;h_1 g \gamma)$ and $f_0 (s_0;h_0 g^{-1} \gamma)$, where $t_1=h_1 s_1$, $t_0=h_0s_0$ and the maps $h_1, h_0\in G$ do not depend on $\psi$ or $\gamma$. Thus the relevant analogue of \eqref{eqLikDeriv} is
			\[
			\Bigl|\frac{\partial t_1}{\partial s_1}\Bigr| \Bigl|\frac{\partial t_0}{\partial s_0}\Bigr| \nabla_\psi L(\psi;\gamma, t_1,t_0)
			\]
			and the sufficient condition \eqref{eqSuffUnbalanced} reduces to
			\[
			\int_{\mathscr{Y}}\int_{\mathscr{Y}} \nabla_\psi\ell(\psi,\lambda\hspace{0.5pt}; t_1, t_0)f_{T_1}(t_1\hspace{0.5pt};\psi, \gamma)f_{T_0}(t_0\hspace{0.5pt};\psi, \gamma)dt_1 dt_0 = 0.
			\]
			Since $T_1$ and $T_0$ have a joint distribution that is parametrized $\psi$-symmetrically with respect to $(\psi,\gamma)$, and such that the log-likelihood function, when expressed in terms of $u_1=g^{-1}t_1$ and $u_0=gt_0$, is antisymmetric, both conditions of Proposition \ref{propConsistency} hold by the proofs of Propositions \ref{propGroup} and \ref{propConsistencyPairs}.
		\end{proof}

		\subsection{Proof of Proposition \ref{propCWPoisson}}
		
		\begin{proof}
			
			Let $S_j=Y_{j0}+Y_{j1}$ and $D_j=Y_{j1}-Y_{j0}$ with realizations $s_j$, $d_j$. Using the notation of Cox and Wong, further define $r_j(\theta)=r_{j0}e^{-\theta}+r_{j1}e^{\theta}$, $r_j(\theta,\nu\omega)=r_j(\theta)+\nu\omega$. The contribution of the $j$th stratum to the log-likelihood function is
			\[
			\log\{\Gamma(s_j+\omega)/\Gamma(\omega)\} + \theta d_j + \omega \log \omega +\omega \log \nu - (s_j+\omega)\log r_j(\theta,\nu\omega)
			\]
			and the log-likelihood function $\ell$ is a sum of $m$ such contributions. On noting that $\mathbb{E}_{\theta,\nu,\omega}(S_j)=r_j(\theta)\mathbb{E}_{\nu,\omega}(\gamma_j)=r_j(\theta)/\nu$ under the assumed model, direct calculation shows that the Fisher information under the gamma random effects assumption has components $\check{\imath}_{\theta\omega}=0=\check{\imath}_{\nu\omega}$ identically in $(\theta,\nu,\omega)$, implying by the clause of Proposition \ref{propCWPoisson} that $\Theta \perp_m \Omega$, $\mathcal{N} \perp_m \Omega$, where $\Theta$, $\Omega$ and $\mathcal{N}$ are the parameter spaces for $\theta$, $\omega$ and $\nu$. The parameter $\omega$ does not therefore play a role, and the problem can be treated as though $\omega$ were known.

			On noting that $\mathbb{E}(D_j)=\nabla_\theta r_j(\theta)/\nu$
			\begin{eqnarray*}
				\mathbb{E}_{m}(\nabla_\theta \ell) = \mathbb{E}_{\theta,\nu,\omega}(\nabla_\theta \ell)&=&   \sum_{j} \Bigl(\frac{\nabla_\theta r_j(\theta)}{\nu} - \frac{\{r_j(\theta)/\nu +\omega\}\nabla_\theta r_j(\theta)}{r_{j}(\theta,\nu\omega)}\Bigr) = 0, \\
				\mathbb{E}_{m}(\nabla_\nu \ell) =\mathbb{E}_{\theta,\nu,\omega}(\nabla_\nu \ell)&=&  \sum_{j} \Bigl( \frac{\omega}{\nu} -\frac{\{r_j(\theta)/\nu +\omega\}\omega}{r_{j}(\theta,\nu\omega)}\Bigr) = 0,
			\end{eqnarray*}
			verifying the weaker condition of Proposition \ref{propCoxWong}.
			
		\end{proof}

		\subsection{Proof of Proposition \ref{propOverstrat}}
		
		\begin{proof}
			
			Index the expectation under the overstratified model by $\xi = (\xi_1,\ldots, \xi_n)=(x_1^\T\psi+w_1^\T\lambda, \ldots, x_n^\T \psi+w_n^\T\lambda)$ with $\xi^* =(x_1^\T\psi^*+w_1^\T\lambda, \ldots, x_n^\T \psi^*+w_n^\T\lambda)$ and true value $\eta = (\eta_1,\ldots, \eta_n)=(x_1^\T\psi^*, \ldots, x_n^\T \psi^*)$. Then, $m$ is synonymous with $\eta$ and we write $\mathbb{E}_{m}=\mathbb{E}_\eta$. With $g_\psi = \mathbb{E}_{\eta}[\nabla_\psi \ell(\psi^*,\lambda)]$ and $g_\lambda = \mathbb{E}_{\eta}[\nabla_\lambda \ell(\psi^*,\lambda)]$,
			\begin{eqnarray*}
				g_\psi &=& \mathbb{E}_{\eta}(S) - \sum_{i=1}^{n}\nabla_\psi K(x_{i}^\T \psi + w_i^\T \lambda) = \sum_{i=1}^n x_i d_i \\
				g_\lambda &=& \mathbb{E}_{\eta}(T) - \sum_{i=1}^{n}\nabla_\lambda K(x_{i}^\T \psi + w_i^\T \lambda) = \sum_{i=1}^n w_i d_i,
			\end{eqnarray*}
			where $d_i=K^\prime (\eta_i) - K^\prime (\xi_i^*)$.
			On letting $D^*=\text{diag}(K^{\prime\prime}(\xi_1^*), \ldots, K^{\prime\prime}(\xi_n^*))$, the components of the information matrix $i$ at $(\psi^*,\lambda)$ are $i_{\psi\psi}=\sum_{i=1}^n K^{\prime\prime}(\xi^*_i)x_ix_i^\T = X^\T D^* X$, $i_{\lambda\lambda}=\sum_{i=1}^n K^{\prime\prime}(\xi^*_i)w_iw_i^\T =W^\T D^* W$, $i_{\psi\lambda}=\sum_{i=1}^n K^{\prime\prime}(\xi^*_i)x_iw_i^\T = X^\T D^* W$ and $i_{\lambda\psi}=\sum_{i=1}^n K^{\prime\prime}(\xi^*_i)w_ix_i^\T = W^\T D^* X$.

			The relevant components of $i^{-1}$ from Proposition \ref{propCoxWong} are
			\begin{eqnarray*}
				i^{\psi\psi} &=& (i_{\psi\psi}-i_{\psi\lambda}i_{\lambda\lambda}^{-1}i_{\lambda\psi})^{-1}  \\
				i^{\psi\lambda}&=&-(i_{\psi\psi}-i_{\psi\lambda}i_{\lambda\lambda}^{-1}i_{\lambda\psi})^{-1}i_{\psi\lambda}i_{\lambda\lambda}^{-1} = i^{\psi\psi} i_{\psi\lambda}i_{\lambda\lambda}^{-1}
			\end{eqnarray*}
			
			If derivatives of $K$ higher than 2 are null, as for linear regression, $d=D^*W\lambda$ by a Taylor expansion of $K^{\prime}(\eta)$ so that $i^{\psi\psi}g_\psi=i^{\psi\psi}X^\T D^* W \lambda$ and $i^{\psi\lambda}g_\lambda = -i^{\psi\psi} X^\T D^* W\lambda$, verifying Proposition \ref{propCoxWong}. More generally however, $d\neq D^*W\lambda$ and the consistency arises instead because $\lambda_m^0=\lambda^*=0$, at which point the condition is trivially satisfied.
		\end{proof}

	\end{appendix}

	\subsection*{Acknowledgements}
	
	This research was partially supported by the Natural Sciences and Engineering Research Council of Canada and the UK Engineering and Physical Sciences Research Council.	We are grateful to the  reviewers for their careful reading and very constructive criticism.

	\bibliographystyle{amsplain}

\begin{thebibliography}{3}
		
		\bibitem {B-N1978}
		\textsc{Barndorff-Nielsen, O.~E.~}(1978).
		\newblock \textit{Information and Exponential Families in Statistical Theory}.
		\newblock Wiley, New York.
		
		\bibitem {B-NC1994}
		\textsc{Barndorff-Nielsen, O.~E.~}and \textsc{Cox,~D.~R.} (1994).
		\newblock \textit{Inference and Asymptotics}.
		\newblock Chapman \& Hall, London.
		
		
		
		\bibitem {Battey2023}
		\textsc{Battey, H.~S.~}(2023).
		Inducement of population sparsity
		\newblock \textit{Canad.~J.~Statist.}, 51, 760--768.
		

		
		
		
		\bibitem {BC2020}
		\textsc{Battey, H.~S.~}and \textsc{Cox, D.~R.~}(2020).
		High dimensional nuisance parameters: an example from parametric survival analysis.
		\newblock \textit{Information Geometry}, 3, 119--148.
		
		\bibitem{B2022}
		\textsc{Battey, H.~S.~}(2022). 
		Discussion of `Assumption-lean inference for generalised linear model parameters' by Vansteelandt and Dukes
		\newblock \emph{J.~Roy.~Statist.~Soc.~Ser.~B}, 84, 696--698.
		
			\bibitem {Breiman2001}
	\textsc{Breiman, L.~}~(2001). Statistical modeling: the two cultures.
	\newblock \emph{Statist.~Sci.}, 16, 199--215.	
		
		\bibitem {Cherno2018}
		\textsc{Chernozhukov, V., Chetverikov, D., Demirer, M.,
			Duflo, E., Hansen, C., Newey, W.~and
			Robins, J.~}(2018). Double/debiased machine learning for treatment and structural parameters
		\newblock \emph{Econom. J.}, 21, 1--68.
		
		
		\bibitem {Cox1961}
		\textsc{Cox, D.~R.}~(1961). Tests of separate families of hypotheses.
		\newblock \emph{Proc. 4th Berkeley Sympos. Vol.~I}, 105--123.
		
		\bibitem {Cox1962}
		\textsc{Cox, D.~R.}~(1962). Further results on tests of separate families of hypotheses
		\newblock \emph{J.~R.~Statist.~Soc.~B}, 24, 406--424.	
		
		\bibitem {Cox2001}
		\textsc{Cox, D.~R.}~(2001). Comment on Breiman's ``Statistical modeling: the two cultures''.
		\newblock \emph{Statist.~Sci.}, 16, 216--218.	
		
		
		
		\bibitem {Cox2012}
		\textsc{Cox, D.~R.}~(2012). Statistical causality: some historical remarks.
		\newblock \emph{in Causality: statistical perspectives and applications}, 1--5.
		
		\bibitem {CoxReid1987}
		\textsc{Cox, D.~R.}~and \textsc{Reid, N.}~(1987). Parameter orthogonality and approximate conditional inference (with discussion).
		\newblock \emph{J.~R.~Statist.~Soc.~B}, 49, 1--39.
		
			\bibitem {CoxWermuth1990}
	\textsc{Cox, D.~R.~}and \textsc{Wermuth, N.}~(1990). 
	An approximation to maximum likelihood estimates in reduced models.
	\newblock \emph{Biometrika}, 77, 747--761.
		
		\bibitem {CoxReid2004}
		\textsc{Cox, D.~R.~}and \textsc{Reid, N.}~(2004). A note on pseudolikelihood constructed from marginal densities.
		\newblock \emph{Biometrika}, 91, 729--737.
		
		\bibitem {CoxWong2010}
		\textsc{Cox, D.~R.~}and \textsc{Wong, M.~Y.}~(2010). A note on the sensitivity of assumptions of a generalized linear mixed model.
		\newblock \emph{Biometrika}, 97, 209--214.
		
		\bibitem {Daniel}
\textsc{Daniel, R., Zhang, J.~}and \textsc{Farewell, D.~}(2021).
Making apples from oranges: comparing noncollapsible effect estimators and their standard errors after adjustment for different covariate sets.
\newblock \textit{Biometrical Journal}, 63, 528--557.
		
		\bibitem {deStavola}
		\textsc{de Stavola, B.~L.~}and \textsc{Cox, D.~R.~}(2008).
		On the consequences of overstratification.
		\newblock \textit{Biometrika}, 95, 992--996.
		
		\bibitem{ED2023}
		\textsc{Evans, R.~J.~}and \textsc{Didelez, V.~}(2023).
		Parameterizing and simulating from causal models (with discussion).
		\newblock \emph{J.~Roy.~Statist.~Soc.~Ser.~B}, to appear.
		
	
		
		\bibitem{Fraser1968}
		\textsc{Fraser, D.~A.~S.~}(1968).
		\newblock \textit{The Structure of Inference}.
		\newblock Wiley, New York.
		
		
	
		
		\bibitem {Huber1967}
		\textsc{Huber, P.~J.}~(1967). The behavior of maximum likelihood estimates under nonstandard conditions.
		\newblock \emph{Proc.~5th Berkeley Sympos. Vol.~I}, 221--233.
		
		\bibitem {Kent1982}
		\textsc{Kent, J.~}~(1982). Robust properties of likelihood ratio tests.
		\newblock \emph{Biometrika}, 69, 19--27.
		
		\bibitem{KW1956}
		\textsc{Kiefer, J.~G.~}and 	\textsc{Wolfowitz, J.~}(1983) Consistency of the maximum likelihood estimator in the presence of many incidental parameters.
		\newblock \textit{Ann.~Math.~Statist.}, 27, 887--906.
		
		
		\bibitem{Lindsay1983}
		\textsc{Lindsay, B.~G.~}(1983)	The geometry of mixture likelihoods: general theory.
		\newblock \textit{Ann.~Statist.}, 11, 86--94.
		
		\bibitem{Lindsay1985}
		\textsc{Lindsay, B.~G.~}(1985).
		Using empirical partially Bayes inference for increased efficiency.
		\newblock \textit{Ann.~Statist.}, 13, 914--931.	
		
		
		\bibitem{Neuhaus1994}
		\textsc{Neuhaus, J.~M.~}and \textsc{Kalbfleisch, J.~D.~}and \textsc{Hauck, W.~W.~}(1994)Conditions for consistent estimation in mexed-effects models for binary matched-pairs data. 
		\newblock \textit{Canad.~J.~Statist.}, 22, 139--148.
		

		
		
		\bibitem {Sch}
		\textsc{Schielzeth, H.~}and \textsc{Dingemanse, N.~}and \textsc{Nakagawa, S.~}and \textsc{Westneat, D.~}and \textsc{Allegue, H.~}and \textsc{Teplitsky, C.~}and \textsc{R{\'e}ale, D.~}and \textsc{Dochtermann, N.~}and \textsc{Garamszegi, L.~}and \textsc{Araya-Ajoy, Y.~}(2020). Robustness of linear mixed-effects models to violations of distributional assumptions.
		\newblock \emph{Methods in Ecology and Evolution}, 11, 1141--1152.
		
			\bibitem{Robins94}
		\textsc{Robins, J.~M.~}, \textsc{Rotnitzky, A.~}and \textsc{Zhao, L.~P.~}(1994). 
		Estimation of regression coefficients when some regressors are not always observed. 
		\newblock \emph{J.~Am.~Stat.~Assoc.}, 89, 846--866.
		
		\bibitem{VD2022}
		\textsc{Vansteelandt, S.~}and \textsc{Dukes, O.}(2022). 
		Assumption-lean inference for generalised linear model parameters (with discussion). 
		\newblock \emph{J.~Roy.~Statist.~Soc.~Ser.~B}, 84, 657--685.
		
		\bibitem {White1982a}
		\textsc{White, H.~}~(1982a). Maximum likelihood estimation of misspecified models.
		\newblock \emph{Econometrica}, 50, 1--25.
		
		
		\bibitem {White1982b}
		\textsc{White, H.~}~(1982b). Regularity conditions for Cox's test of non-nested hypotheses.
		\newblock \emph{J.~Econometrics}, 19, 301--318.
		

		
		
		
		
		
	\end{thebibliography}

	\end{document}